\documentclass[twocolumn]{IEEEtran}

\usepackage{color}
\usepackage{paralist, comment}
\usepackage{booktabs, multirow}
\usepackage{pstool}
\usepackage{amsmath,amsfonts,amssymb}
\usepackage{mathabx,fixmath}
\usepackage{algorithm}
\usepackage{algorithmic}
\usepackage{cite}
\usepackage[lowtilde]{url}
\def\x{\mathbold{x}}
\def\A{\mathbold{A}}

\def\Q{\mathbf{Q}}

\def\c{\mathbold{c}}

\def\y{\mathbold{y}}

\newcommand{\m}{\mathbold{m}}

\newcommand{\p}{\mathbold{p}}

\def\X{\mathbf{X}}

\def\transp{\mathsf{T}}

\newcommand{\vv}{\mathbold{\mu}}
\def\qed{\hfill $\blacksquare$ }
\DeclareMathOperator*{\minimize}{\textrm{minimize}}
\newtheorem{remark}{Remark}

\newtheorem{proposition}{Proposition}
\newtheorem{assumption}{Assumption}

\newtheorem{theorem}{Theorem}

\input{mysymbol.sty}

\newcommand{\bm}{{\bf m}}

\title{Prediction-Correction Algorithms for\\ Time-Varying Constrained Optimization}

\author{Andrea~Simonetto and Emiliano~Dall'Anese 
\thanks{This paper significantly expands the results and presents convergence proofs that are referenced in \emph{Proceedings of the Twentieth IFAC World Congress}, July 2017, \cite{Simonetto2017}. 

A. Simonetto is with IBM Research, Ireland, Dublin, Ireland. E. Dall'Anese is  with the National Renewable Energy Laboratory, Golden, CO, USA.  E-mail: {andrea.simonetto@ibm.com}, {emiliano.dallanese@nrel.gov}.%
}
}

\begin{document}

\maketitle

\begin{abstract}%
This paper develops online algorithms to track solutions of time-varying constrained optimization problems. Particularly, resembling workhorse Kalman filtering-based approaches for dynamical systems, the proposed methods involve prediction-correction steps to provably track the trajectory of the optimal solutions of time-varying convex problems. The merits of existing prediction-correction methods have been shown for unconstrained problems and for setups where computing the inverse of the Hessian of the cost function is computationally affordable. This paper addresses the limitations of existing methods by tackling constrained problems and by designing first-order  prediction steps that rely on the Hessian of the cost function (and do not require the computation of its inverse). In addition, the proposed methods are shown to improve the convergence speed of existing prediction-correction methods when applied to unconstrained problems.  Numerical simulations corroborate the analytical results and showcase performance and benefits of the proposed algorithms. A realistic application of the proposed method to real-time control of energy resources is presented.

\end{abstract}%



\section{Introduction}





We consider constrained optimization problems that vary continuously in time. We outline the main ideas by first focusing on problems with a time-varying objective function. Consider then the following optimization problem:
\begin{equation}\label{eq.problem}
\x^*(t) :=\argmin_{\x\in X} f(\x; t), \quad \textrm{for } t \geq 0 \;, 
\end{equation}
where $X \subseteq \reals^n$ is a convex set; $t\in \reals_+$ is non-negative, continuous, and it is used to index time; and,  $f: \mathbb{R}^n\times \mathbb{R}_{+} \to \mathbb{R}$ is a \emph{smooth strongly convex function}.
The goal is to find (and track) the solution $\x^*(t)$ of~\eqref{eq.problem} for each time $t$ -- hereafter referred to as the  optimal solution \emph{trajectory}. 

Problem~\eqref{eq.problem} is a generalization of traditional time-invariant (i.e., static) problems, and can naturally model: a)  control problems where one seeks to generate a control action depending on a (parametric) varying optimization problem~\cite{Jerez2014,Hours2014,Gutjahr2016}; b) signal processing problems~\cite{Jakubiec2013}, where states of a dynamical process are  estimated online based on time-varying observations -- including time-varying compressive sensing settings~\cite{Asif2014,Yang2015,Vaswani2015,Balavoine2015,Simonetto2015a,Sopasakis2016}; and, c)  inferential problems on dynamic networks~\cite{Baingana2015}. Additional application domains include  robotics\cite{ardeshiri2010convex,verscheure2009time,Koppel2015a}, smart grids~\cite{Zhao2014,DallAnese2016,DallAnese2016a}, economics~\cite{Dontchev2013}, real-time magnetic resonance imaging (MRI)~\cite{Uecker2012}, and human migration studies~\cite{Nagurney2006}.

The problem \eqref{eq.problem} might be solved in a centralized setting based on a continuous time platform \cite{Ye2015,Rahili2015a,Rahili2015,Gong2016,Fazlyab2015,Fazlyab2016}; however, here we focus on a discrete time setting. The reason for this choice is motivated by the widespread use of digital computing units, such as control units (actuators) and digital sensors. In this context, we envision that our optimization problem will change in response to measurements taken at discrete time steps and its solution could provide control actions to be implemented on digital control units. 

Therefore, we use sampling arguments to reinterpret~\eqref{eq.problem}  as a sequence of time-invariant problems. In particular, upon sampling the objective functions $ f(\x; t)$ at time instants $t_{k}$,  $k=0,1,2,\dots$, where the sampling period $h := t_k - t_{k-1}$ can be chosen arbitrarily small, one can solve the sequence of time-invariant problems
\begin{equation}\label{eq.problemd}
\x^*(t_k) :=\argmin_{\x\in X} f(\x; t_k) ,  \hspace{.3cm}  k \in \mathbb{N}. 
\end{equation}
By decreasing $h$, an arbitrary accuracy may be achieved when approximating problem~\eqref{eq.problem} with \eqref{eq.problemd}. However, solving \eqref{eq.problemd} for each sampling time $t_k$ may not be computationally affordable in many application domains, even for moderate-size problems\footnote{Note that here we focus on constant sampling periods; variable ones induced by self-triggering have been considered in unconstrained problems~\cite{Fazlyab2016a} and could be considered in constrained problems in future works.}. 
%
%


Focusing on unconstrained optimization problems, the works in~\cite{Paper1, Paper2} developed a prediction-correction method to find and track the solution trajectory $\x^*(t)$ up to a bounded asymptotical error, starting from an arbitrary guess $\x_0$. This methodology arises from non-stationary optimization~\cite{Polyak1987,Popkov2005}, parametric programming~\cite{Robinson1980,Dontchev2009,Zavala2010,Dontchev2013}, and continuation methods in numerical mathematics~\cite{Allgower1990}. It also resembles evolutionary variational inequalities~\cite{Cojocaru2005,Nagurney2006} and path-following methods in interior point solvers~\cite{Nesterov2012}. This paper significantly broadens the method of~\cite{Paper1} by offering the following contributions. 
\begin{enumerate}
\item[\emph{(i)}] We develop prediction-correction methods to track the solutions of the time-varying \emph{constrained} problems~\eqref{eq.problem}; 
\item[\emph{(ii)}] We develop first-order algorithms that do not involve the computation of the inverse of the Hessian of the cost function, as required in~\cite{Paper1}; the proposed prediction-correction method is computationally lighter, as it requires only matrix-vector multiplications. Further,  we offer a trade-off between tracking capabilities  and computational effort; and
\item[\emph{(iii)}] We show that the proposed prediction-correction algorithm improves on the method in~\cite{Paper1} when applied to unconstrained optimization problems; particularly, it exhibits enhanced local convergence properties by relying on Newton-like prediction steps. 
\end{enumerate}
The design and analysis of proposed prediction-correction methods are grounded on the theory of generalized equations and implicit function theorems~\cite{Dontchev2009}. 


{\bf Organization.} In Section~\ref{sec:probform}, we describe the prediction-correction methodology for constrained time-varying optimization problems. The special case of unconstrained case is developed in Section~\ref{sec:unc}. Convergence analysis is discussed in Section~\ref{sec:conv}, both for global $O(h)$ convergence and local $O(h^2)$ convergence. Numerical examples are displayed in Section~\ref{sec:num} -- including a realistic application in real-time control of energy resources --, and in Section~\ref{sec:concl}, we draw our conclusions. The proofs of all the propositions and theorems are given in the appendices.


{\bf Notation.} Vectors are written as $\x\in\reals^n$ and matrices as $\A\in\reals^{n\times n}$. We use $\|\cdot\|$ to denote the Euclidean norm in the vector space, and the respective induced norms for matrices and tensors. The gradient of the function $f(\x; t)$ with respect to $\x$ at the point $(\x,t)$ is denoted as $\nabla_{\x} f(\x; t) \in \reals^n$, while the partial derivative of the same function with respect to (w.r.t.) $t$ at $(\x,t)$ is written as $\nabla_t f(\x; t)\in \reals$. Similarly, the notation $\nabla_{\x\x} f(\x; t) \in \reals^{n\times n}$ denotes the Hessian of $f(\x;t)$ w.r.t. $\x$ at $(\x,t)$, whereas $\nabla_{t\x} f(\x; t) \in \reals^{n}$ denotes the partial derivative of the gradient of $f(\x;t)$ w.r.t. the time $t$ at $(\x,t)$, i.e. the mixed first-order partial derivative vector of the objective. The tensor $\nabla_{\x\x\x} f(\x; t) \in \reals^{n\times n\times n}$ indicates the third derivative of $f(\x;t)$ w.r.t. $\x$ at $(\x,t)$, the matrix $\nabla_{\x t\x} f(\x; t) = \nabla_{t\x\x} f(\x; t)   \in \reals^{n\times n}$ indicates the time derivative of the Hessian of $f(\x;t)$ w.r.t. the time $t$ at $(\x,t)$, and the vector $\nabla_{t t\x} f(\x; t) \in \reals^{n}$ indicates the second derivative in time of the gradient of $f(\x;t)$ w.r.t. the time $t$ at $(\x,t)$. The indicator function is indicated as $\iota_{X}(\x)$ for the convex set $X \subseteq \mathbb{R}^n$; by definition $\iota_{X}(\x) = 0$ when $\x \in X$ and $+\infty$ otherwise. The subdifferential of $\iota_{\X}$ is the set-valued map known as normal cone $N_{X}: \mathbb{R}^n \rightrightarrows \mathbb{R}^n$ (Cf.~\cite{Dontchev2009}).


\section{Prediction-Correction Strategy}\label{sec:probform}

In this section, we first focus on problem~\eqref{eq.problem} and design a prediction-correction algorithm to track the (unique) trajectory of the optimal solution. As explained in the introduction, consider sampling~\eqref{eq.problem} at times $t_k$, $k \in \mathbb{N}$, and constructing a sequence of time-invariant problems~\eqref{eq.problemd}. In lieu of solving ~\eqref{eq.problemd} at each time step, the goal of the prediction-correction strategy is to determine an approximate optimizer for~\eqref{eq.problem} at $t_{k+1}$ in a computationally affordable way, given the current approximate optimizer at $t_k$. 

\subsection{Prediction}

Suppose that $\x_k$ is an approximate solution of~\eqref{eq.problemd} at time $t_k$. Given $\x_k$,  the prediction step seeks an approximate optimizer for~\eqref{eq.problem} at $t_{k+1}$, given  the only information available at time $t_{k}$. Let $\x_{k+1|k}$ denote the output of the prediction step, which is computed as explained next. 

Notice first that solving the time-invariant problem~\eqref{eq.problemd} associated with time $t_k$ is equivalent to solving the generalized equation
\begin{equation}\label{ge}
\nabla_{\x} f(\x^*(t_k); t_k) + N_{X}(\x^*(t_k)) \ni \mathbf{0}
\end{equation}
where $N_X$ is the normal cone operator, while $\x^*(t_k)$ is the optimizer of~\eqref{eq.problemd} at $t_k$. Although many equivalent ways to formulate the solution of~\eqref{eq.problemd} exist, e.g., as the fixed-point of a properly defined operator, here we focus on generalized equations to build on the implicit function theorems of~\cite{Dontchev2013,Dontchev2009}.  

With~\eqref{ge} in place, the prediction step seeks the solution of the following perturbed generalized equation
\begin{multline}\label{ger0}
\nabla_{\x} f(\x_{k+1|k}; t_{k+1}) + N_{X}(\x_{k+1|k}) \approx \\
\nabla_{\x} f(\x_{k}; t_{k}) + \nabla_{\x\x} f(\x_{k}; t_{k}) (\x_{k+1|k} - \x_k)  \\ +  h\, \nabla_{t\x} f(\x_{k}; t_{k}) + N_{X}(\x_{k+1|k}) \ni \mathbf{0}.
\end{multline}
That is, the prediction step produces a solution that is optimal w.r.t. a perturbed (first-order) version of the original generalized equation~\eqref{ge}. The idea of this perturbed approximate version is motivated by the fact that one would like to solve $\nabla_{\x} f(\x_{k+1|k}; t_{k+1}) + N_{X}(\x_{k+1|k})\ni \mathbf{0}$, but that is not possible at time $t_{k}$. Therefore, we perturb said generalized equation by a backward Taylor expansion, in a way that is possible to solve it at time $t_k$.

We can now replace~\eqref{ger0} with the following equivalent formulation
\begin{multline}\label{qp}
\x_{k+1|k} \!= \!\argmin_{\x \in X} \Big\{\frac{1}{2} \x^\transp \nabla_{\x\x} f(\x_k; t_k) \x \\ + (\nabla_{\x} f(\x_k; t_k) + h\, \nabla_{t\x} f(\x_k; t_k) -\nabla_{\x\x} f(\x_k; t_k)\x_k )^\transp\x\Big\},
\end{multline}
where we can notice that the normal cone operator yields the feasible set over which the optimization is computed. 

In lieu of seeking an exact solution of~\eqref{qp} (which is a constrained optimization problem with quadratic cost), consider the less computational demanding task of finding an approximate solution of~\eqref{qp} by computing a number of  projected gradient descent steps -- the first key step towards a \emph{first-order} prediction-correction method (Note that we consciously avoid Newton's method here for the high computational complexity burden of computing the Hessian inverse). Particularly, let ${\widehat{\x}^0}$ be a dummy variable initialized as ${\widehat{\x}^0} = \x_k$; then, the following steps are performed:
\begin{multline}\label{gpred}
{\widehat{\x}^{p+1}} = \mathbb{P}_{X}[\widehat{\x}^{p} -  \alpha (\nabla_{\x\x} f(\x_k; t_k)(\widehat{\x}^{p}-\x_k) \\+ h\, \nabla_{t\x} f(\x_k; t_k) + \nabla_{\x} f(\x_k; t_k))],
\end{multline}
for $p = 0, 1, \dots, P-1$, where $P$ is a pre-determined number of gradient steps, $\alpha > 0$ is the stepsize, and $\mathbb{P}_X$ is the projection operator over the convex set $X$. Once $P$ steps are performed, $\tilde{\x}_{k+1|k}$ is set to:
\begin{equation}
\tilde{\x}_{k+1|k} = \widehat{\x}^{P}.
\end{equation}

\subsection{Correction}

Once the cost function $f(\cdot; t_{k+1})$ becomes available, the correction step is performed to refine the estimate of the optimal solution $\x^*(t_{k+1})$. To this end, a first-order projected gradient method is considered next. Particularly, let $\widehat{\x}^{0} = \tilde{\x}_{k+1|k}$ be  a dummy variable; then, consider the following projected gradient steps
\begin{equation}\label{gcorr}
{\widehat{\x}^{c+1}} \!\!= \mathbb{P}_{X}[{\widehat{\x}^{c}} - \\ \beta (\nabla_{\x} f(\widehat{\x}^{c}; t_{k+1})],
\end{equation}
for $c = 0, 1, \dots, C-1$, where $C$ is a predetermined number of gradient steps and $\beta >0 $ the stepsize. The estimate of the optimal solution $\x^*_{k+1}$ is then computed as 
${\x}_{k+1} = \widehat{\x}^{C}$.


\subsection{Complete Algorithm}

The complete algorithm Constrained - First Order Prediction Correction (C-FOPC) is tabulated as Algorithm~\ref{algo_fo}. Steps 4-7 are utilized to compute $ \tilde{\x}_{k+1|k}$ based on the information available at $t_{k}$. Provided that the projection operator is easy to carry out (set $X$ is simple, e.g., it is an $\ell_2$ or $\ell_{\infty}$ norm), and the Hessian is easy to evaluate, the computational complexity of these steps is $O(Pn^2)$, which is quadratic (due to matrix-vector multiplications) in the number of scalar decision variables. This is in contrast with  the algorithms presented in \cite{Paper1}, which involve the computation of the Hessian inverse and whose computational complexity is therefore $O(n^3)$. Steps 10-14 are utilized to compute ${\x}_{k+1}$, based on the information available at $t_{k+1}$. Provided that the projection operator is easy to perform (set $X$ is simple) and the gradient is easy to evaluate, the computational complexity of these steps is $O(Cn)$, which is linear in the number of scalar decision variables. 


\begin{remark}
\emph{[Time derivative approximation]} The time derivative of the gradient $\nabla_{t\x} f(\x; t)$ can be substituted with an approximate version, as explained in~\cite{Paper1}.
\end{remark}

{\begin{algorithm}[tb]
\caption{Constrained First-Order Prediction-Correction (C-FOPC)}\label{algo_fo} 
\begin{algorithmic}[1] 
\small{\REQUIRE  Initial variable $\x_0$. Initial objective function $f(\x;t_0)$, no. of prediction steps $P$ and correction steps $C$
\FOR {$k=0,1,2,\ldots$}
\STATE \color{blue}// time $t_k$ \color{black}
   \STATE Prediction: initialize ${\widehat{\x}^0} = \x_k$
   \FOR {$p=0:P-1$}
          \STATE  Predict the variable by the  gradient step  [cf \eqref{gpred}]   
         \begin{multline*}
{\widehat{\x}^{p+1}} \!\!= \mathbb{P}_{X}[ {\widehat{\x}^{p}} -  \alpha (\nabla_{\x\x} f(\x_k; t_k){(\widehat{\x}^{p} - \x_k)} +\\ h\, \nabla_{t\x} f(\x_k; t_k) + \nabla_{\x}f(\x_k; t_k))] \end{multline*}
\ENDFOR
\STATE Set the predicted variable $\tilde{\x}_{k+1|k}=\widehat{\x}^{P}$
\STATE \color{blue}// time $t_{k+1}$ \color{black}
\STATE Acquire the updated function $f(\x;t_{k+1})$
   \STATE Initialize the sequence of corrected variables $\widehat{\x}^{0} = \tilde{\x}_{k+1|k}$ 
     \FOR {$c=0:C-1$}
          \STATE  Correct the variable by the  gradient step  [cf \eqref{gcorr}]   
         $${\widehat{\x}^{c+1}} \!\!= \mathbb{P}_{X}[{\widehat{\x}^{c}} - \beta (\nabla_{\x} f(\widehat{\x}^{c}; t_{k+1})]$$
          \ENDFOR
   \STATE Set the corrected variable ${\x}_{k+1}=\widehat{\x}^{C}$
\ENDFOR}
\end{algorithmic}\end{algorithm}
%

\subsection{Special Case: Unconstrained Problems}\label{sec:unc}

In this section, we focus on the special case of unconstrained problems (i.e., $X = \reals^n$). Although~\cite{Paper1} has given an extensive characterization of methods for unconstrained problems, we will see in this section that further important improvements can be achieved based on the prediction generalized equation~\eqref{ger0}. 

For unconstrained problems, the prediction equation~\eqref{ger0} can be rewritten as
\begin{multline}\label{ger0.0}
\nabla_{\x} f(\x_{k+1|k}; t_{k+1}) \approx 
\nabla_{\x} f(\x_{k}; t_{k}) + \\ \nabla_{\x\x} f(\x_{k}; t_{k}) (\x_{k+1|k} - \x_k)  +  h\, \nabla_{t\x} f(\x_{k}; t_{k}) = \mathbf{0}.
\end{multline}
Since the suboptimality at time $t_k$ can be easily characterized by the gradient, i.e., $\nabla_{\x} f(\x_{k}; t_{k})$, we can also modify~\eqref{ger0.0} to tune the prediction step. A way to do that is to require $\x_{k+1|k}$ to reduce the suboptimality by a factor of $1-\gamma \in [0,1]$, and pose the problem as 
\begin{multline}\label{ger0.1}
\nabla_{\x} f(\x_{k+1|k}; t_{k+1}) \approx 
\nabla_{\x} f(\x_{k}; t_{k}) + \nabla_{\x\x} f(\x_{k}; t_{k})\times\\ (\x_{k+1|k} - \x_k)  +  h\, \nabla_{t\x} f(\x_{k}; t_{k}) = (1-\gamma)\nabla_{\x} f(\x_k; t_{k}).
\end{multline}
When $\gamma = 1$,~\eqref{ger0.1} boils down to~\eqref{ger0.0} (when one seeks an ``optimal prediction''); on the other hand, when $\gamma = 0$, ~\eqref{ger0.1} coincide with~\cite{Paper1}, where a prediction vector that maintains the same suboptimality (i.e., the same gradient) between successive time steps is sought. Notice that the possibility of tuning the algorithm via $\gamma$ is possible only for unconstrained optimization problems, since we have access to a ``suboptimality measure''. 

From~\eqref{ger0.1} we obtain the prediction vector via the following update
\begin{multline}
\label{pred_up}
\x_{k+1|k} = \x_{k} - \nabla_{\x\x} f(\x_{k}; t_{k})^{-1}\times\\(h\, \nabla_{t\x} f(\x_{k}; t_{k}) + \gamma \, \nabla_{\x} f(\x_{k}; t_{k})),
\end{multline}
which combines a Newton-like step in the direction on the changing cost function, and a (damped) Newton's step towards the optimizer at time $t_k$. When $\gamma = 0$,~\eqref{pred_up} coincides with the prediction step in~\cite{Paper1}. To obtain a first-order update, it is then easy to modify~\eqref{gpred} as
\begin{multline}\label{gpred.0}
{\widehat{\x}^{p+1}} \!\!= \widehat{\x}^{p} - \alpha (\nabla_{\x\x} f(\x_k; t_k)(\widehat{\x}^{p}-\x_k) \\ + h\, \nabla_{t\x} f(\x_k; t_k)+\gamma\,\nabla_{\x} f(\x_k; t_k)),
\end{multline}
where the notation is the same as in~\eqref{gpred}. On the other hand, the correction update becomes in this case
\begin{equation}\label{gcorr.0}
{\widehat{\x}^{c+1}} \!\!= {\widehat{\x}^{c}} - \\ \beta (\nabla_{\x} f(\widehat{\x}^{c}; t_{k+1}).
\end{equation}
The resultant Unconstrained - First Order Prediction Correction (U-FOPC) method  is tabulated as Algorithm~\ref{algo_fou}. 

{\begin{algorithm}[tb]
\caption{Unconstrained First-Order Prediction-Correction (U-FOPC)}\label{algo_fou} 
\begin{algorithmic}[1] 
\small{\REQUIRE  Initial variable $\x_0$. Initial objective function $f(\x;t_0)$, no. of prediction steps $P+1$ and correction steps $C+1$, sub-optimality requirement $\gamma\in[0,1]$
\FOR {$k=0,1,2,\ldots$}
\STATE \color{blue}// time $t_k$ \color{black}
   \STATE Prediction: initialize ${\widehat{\x}^0} = \x_k$
   \FOR {$p=0:P-1$}
          \STATE  Predict the variable by the  gradient step  [cf \eqref{gpred}]   
         \begin{multline*}
{\widehat{\x}^{p+1}}  = {\widehat{\x}^{p}} -\alpha (\nabla_{\x\x} f(\x_k; t_k){(\widehat{\x}^{p} - \x_k)} \\+ h\, \nabla_{t\x} f(\x_k; t_k)+\gamma\,\nabla_{\x} f(\x_k; t_k)) \end{multline*}
\ENDFOR
\STATE Set the predicted variable $\tilde{\x}_{k+1|k}=\widehat{\x}^{P}$
\STATE \color{blue}// time $t_{k+1}$ \color{black}
\STATE Acquire the updated function $f(\x;t_{k+1})$
   \STATE Initialize the sequence of corrected variables $\widehat{\x}^{0} = \tilde{\x}_{k+1|k}$ 
     \FOR {$c=0:C-1$}
          \STATE  Correct the variable by the  gradient step  [cf \eqref{gcorr}]   
         $${\widehat{\x}^{c+1}} \!\!= \mathbb{P}_{X}[{\widehat{\x}^{c}} - \beta (\nabla_{\x} f(\widehat{\x}^{c}; t_{k+1})]$$
          \ENDFOR
   \STATE Set the corrected variable ${\x}_{k+1}=\widehat{\x}^{C}$
\ENDFOR}
\end{algorithmic}\end{algorithm}


\section{Convergence Analysis}\label{sec:conv}

In this section, we establish analytical results to bound the discrepancy between the optimal solution $\x^*(t)$ and the iterates $\x_k$ produced by the prediction-correction schemes developed in Section \ref{sec:probform}. Particularly, we will show that $\x_k$ tracks $\x^*(t)$ up to an error term that depends on the discrete-time sampling period. To this end, some technical conditions are required as stated next


\vskip2mm

\begin{assumption} \label{as.str} {\bf
The function $f(\x; t)$ is twice differentiable and $m$-strongly convex in $\x\in X$ and uniformly in $t$; that is, the Hessian of $f(\x; t)$ with respect to $\x$ is bounded below by $m$ for each $\x\in X$ and uniformly in $t$, 
$$
\nabla_{\x\x} f(\x; t) \succeq m \mathbf{I}, \quad \forall \x\in X, t.
$$  
In addition, the function $f(\x; t)$ has bounded second and third order derivatives with respect to $\x\in X$ and $t$:
$$
\|\nabla_{\x\x} f(\x; t)\|\leq L, \, \|\nabla_{t\x} f(\x; t)\|\leq C_0.
$$
}
\end{assumption}
%

\begin{assumption} \label{as.smooth} 
{\bf
The function $f(\x; t)$ has bounded third order derivatives with respect to $\x\in X$ and $t$:
$$
\|\nabla_{\x\x\x} f(x; t)\|\!\leq \!C_1, \|\nabla_{\x t\x} f(\x; t)\|\!\leq\! C_2,  \|\nabla_{t t\x} f(\x; t)\|\!\leq \!C_3.
$$ 
}
\end{assumption}
\vskip2mm

%

Assumption~\ref{as.str} guarantees that problem~\eqref{eq.problem} is strongly convex and has a \emph{unique} solution for each time instance. Uniqueness of the solution implies that the solution trajectory is also unique. Assumption\ref{as.str} also ensures that the Hessian is bounded from above; this property  is equivalent to the Lipschitz continuity of the gradient. This setting is common in the the time-varying optimization domain; see, for instance~\cite{Popkov2005, Dontchev2013, Jakubiec2013, Ling2013, Paper1, Paper2}. Assumption~\ref{as.smooth} ensures that the third derivative tensor $\nabla_{\x\x\x} f(x; t)$ is  bounded above (typically required for the analysis of Newton-type algorithms). Assumptions~\ref{as.str}-\ref{as.smooth} impose boundedness of the temporal variability of gradient and Hessian. These last properties ensure the possibility to build a prediction scheme based on the knowledge of (or an estimate of) how the function and its derivatives change over time. A similar assumption was required (albeit only locally) for the local convergence analysis in~\cite[Eq.~(3.2)]{Dontchev2013}. 

Assumption~\ref{as.str} is sufficient to show that the solution \emph{mapping} $t \mapsto \x^*(t)$ is single-valued and locally Lipschitz continuous in $t$; in particular, from \cite[Theorem 2F.10]{Dontchev2009} we have that
\begin{equation}\label{eq.lip}
\|\x^*(t_{k+1}) \!-\! \x^*(t_{k})\| \leq \! \frac{1}{m}\|\nabla_{t\x} f(\x; t)\| (t_{k+1}\!-\!t_{k}) \leq \frac{C_0 h }{m}, 
\end{equation} 
for sufficiently small sampling periods $h$. This result  established link between the sampling period $h$ and the temporal variability of the optimal solutions; further,~\eqref{eq.lip} will be utilized to substantiate convergence and tracking  capabilities of the proposed prediction-correction methods.

\subsection{Constrained algorithm global $O(h)$ convergence}

We study the convergence properties of the sequence $\{\x_k\}_{k \in \mathbb{N}}$ generated by the algorithm C-FOPC, for different choices of the stepsize. In the following theorem,  we show that the optimality gap $\|\x_k-\x^*(t_{k})\|$ converges exponentially to a given error bound, globally, i.e., from an arbitrary $\x_0$. 

\vskip2mm
\begin{theorem}\label{theo.constrained_0}

{\bf
Consider the sequence $\{\x_k\}_{k \in \mathbb{N}}$ generated by the C-FOPC method, and let Assumption~\ref{as.str} hold true. Let $\x^*(t_k)$ be the optimizer of~\eqref{eq.problem} at time $t_k$. Define the following quantities 
\begin{equation}\label{stepsize_0}
\varrho_{\mathrm{P}} = \max\{|1-\alpha m|,|1- \alpha L|\}, \,\varrho_{\mathrm{C}} = \max\{|1-\beta m|,|1- \beta L|\},
\end{equation}
and choose the stepsizes $\alpha$ and $\beta$ as
\begin{equation}
\alpha < 2/L, \quad \beta < 2/L.
\end{equation}
Further, select the number of correction steps $C$ in a way that  
\begin{equation}\label{eq.cvrt}
\tau_0 := \varrho_{\mathrm{C}}^{C}\left[\varrho_{\mathrm{P}}^{P} + (\varrho_{\mathrm{P}}^{P}+1) \frac{2L}{m}\right]  <1.
\end{equation} 
Then, the sequence $\{\|\x_k-\x^*(t_k)\|\}_{k \in \mathbb{N}}$ converges linearly with rate $\tau_0$ to an asymptotical error bound, and 
\begin{equation}\label{main.result1_0}
\limsup_{k\to\infty}\|\x_{k} - \x^*(t_k)\| = O(\varrho_{\mathrm{C}}^{C} h).
\end{equation}
}
\end{theorem}

\vskip2mm
\begin{IEEEproof}
See Appendix~\ref{ap.grad}, where we also derive all the constants in the right-hand side of~\eqref{main.result1_0}.
\end{IEEEproof}

\vskip2mm

Theorem~\ref{theo.constrained_0} says that the sequence $\{\x_k\}_{k \in \mathbb{N}}$ generated by C-FOPC converges to a neighborhood of the optimal solution trajectory $\x^*(t)$. In particular, for a choice of prediction and correction steps $P$ and $C$, the sequence converge linearly to the asymptotic bound. The bound depends linearly on the sampling period and exponentially on the correction steps $C$. When one performs the correction step exactly, i.e., $C \to \infty$, then the asymptotic bound goes to zero (in fact, in that case each of the time-invariant problems is solved exactly). Finally, if one looks at the condition~\eqref{eq.cvrt}, one discovers that for global convergence one needs more that standard $\varrho_{\mathrm{C}}^{C} < 1$: in particular $C$ has to be chosen high enough to counteract any errors that come from the prediction step. We see next that this condition can be weakened in the case in which the cost function has higher order smoothness properties like the ones required by Assumption~\ref{as.smooth}.


\subsection{Constrained algorithm local $O(h^2)$ convergence}

We study now the local convergence properties of the sequence $\{\x_k\}_{k \in \mathbb{N}}$ generated by the algorithm C-FOPC under higher order smoothness (Assumption~\ref{as.smooth}), and we show that starting from a small enough optimality gap $\|\x_k-\x^*(t_{k})\|$, and for a small enough sampling period $h$, the error bound can be reduced from $O(h)$ to $O(h^2)$.

\vskip2mm

\begin{theorem}\label{theo.constrained}
{\bf
Consider the sequence $\{\x_k\}_{k \in \mathbb{N}}$ generated by the C-FOPC method, and let Assumptions~\ref{as.str}-\ref{as.smooth} hold true. Let $\x^*(t_k)$ be the optimizer of~\eqref{eq.problem} at time $t_k$. Define the following quantities 
\begin{equation}\label{stepsize}
\varrho_{\mathrm{P}} = \max\{|1-\alpha m|,|1- \alpha L|\}, \, \varrho_{\mathrm{C}} = \max\{|1-\beta m|,|1- \beta L|\},
\end{equation}
and choose the stepsizes $\alpha$ and $\beta$ as
\begin{equation}
\alpha < 2/L, \quad \beta < 2/L.
\end{equation}
Further, select $\tau\in (0,1)$, the number of prediction steps $P$, and the number of correction steps $C$ in a way that  $\varrho_{\mathrm{P}}^{P}\varrho_{\mathrm{C}}^{C} < \tau$. 

There exist an upper bound on the sampling period $\bar{h}$ and a convergence region $\bar{R}$, such that if the sampling period is chosen as $h \leq \bar{h}$ and the initial optimality gap satisfy $\|\x_{0} - \x^*(t_0)\| \leq \bar{R}$, then the sequence $\{\|\x_k-\x^*(t_k)\|\}_{k \in \mathbb{N}}$ converges linearly with rate $\tau$ to an asymptotical error bound, and 
\begin{equation}\label{main.result1}
\limsup_{k\to\infty}\|\x_{k} - \x^*(t_k)\| = O(h^2\,\varrho_{\mathrm{C}}^{C}) + O(h\,\varrho_{\mathrm{P}}^{P}\varrho_{\mathrm{C}}^{C}).
\end{equation}
In addition, the bounds $\bar{h}$ and $\bar{R}$ are given as
\begin{equation}
\bar{h} \!=\! \frac{\tau - \varrho_{\mathrm{C}}^{C}\varrho_{\mathrm{P}}^{P}}{\varrho_{\mathrm{C}}^{C}(\varrho_{\mathrm{P}}^{P} + 1)}\Big(\frac{C_1 C_0}{m^2} + \frac{C_2}{m} \Big)^{-1}\!\!\!, \bar{R} \!=\! \frac{2 \,m}{C_1}\! \Big(\frac{C_1 C_0}{m^2} + \frac{C_2}{m}\Big)\! (\bar{h}-h).
\end{equation}
}
\end{theorem}

\vskip2mm
\begin{IEEEproof}
See Appendix~\ref{ap.grad}, where we also derive all the constants in the right-hand side of~\eqref{main.result1}.
\end{IEEEproof}

\vskip2mm
Theorem~\ref{theo.constrained} asserts that the sequence $\{\x_k\}_{k \in \mathbb{N}}$ generated by C-FOPC locks to a neighborhood of the optimal solution trajectory $\x^*(t)$. In particular, for a choice of prediction and correction steps $P$ and $C$, there exist an upper bound on the sampling period and an attraction region, such that if the sampling period is smaller than the bound and the initial optimality gap is in the attraction region, then the sequence converge (at least) linearly to an asymptotic bound. The bound depends on the sampling period and on the selection of prediction and correction steps $P$ and $C$. When one performs an optimal prediction (that is $P \to \infty$), then the bound goes as $O(h^2)$, which is similar to the bounds derived in~\cite{Paper1}. When one performs the correction step exactly, i.e., $C \to \infty$, then the asymptotic bound goes to zero (in fact, in that case each time-invariant problem is solved exactly). 


The presence of an attraction region is due to mimicking a Newton's step in the prediction stage (that is, the presence of the gradient $\nabla_{\x}f(\x_k; t_k)$ in the generalized equation~\eqref{ger0}). To further understand this point, consider the unconstrained case: as expressed in Section~\ref{sec:unc}, the prediction step optimization problem~\eqref{qp} has a closed-form solution as Eq.~\eqref{pred_up}. The presence of the gradient $\nabla_{\x}f(\x_k; t_k)$ translates into the presence of the Hessian inverse in~\eqref{pred_up}, which is a Newton's step. Our prediction strategy, even in the constrained case, mimics such step and inherits its pros and cons, e.g., the presence of a local convergence region. In addition, as in a Newton's method, when the function is quadratic, then $C_1 = 0$, and the convergence is global. 

By putting together the results of Theorems~\ref{theo.constrained_0}-\ref{theo.constrained}, one can arrive at simple strategies for parameter selection for the designer. In particular, for convergence, the designer needs only to know the strong convexity constant $m$ and the Lipschitz constant $L$, which is a standard requirement. Then, they can design the stepsizes $\alpha$ and $\beta$ as in~\eqref{stepsize_0} (or equivalently~\eqref{stepsize}) and the number of prediction and correction steps as in~\eqref{eq.cvrt}. With this in place, in the worst case they obtain a $O(\varrho_{\textrm{C}}^C h)$ error bound, and in the best case -- if the function satisfies Assumption~\ref{as.smooth}, the optimality gap enters in the attraction region, and $h$ is small enough -- they obtain a better $O(h^2\,\varrho_{\mathrm{C}}^{C}) + O(h\,\varrho_{\mathrm{P}}^{P}\varrho_{\mathrm{C}}^{C})$ error bound. This is the case, since the conditions on $C$ and $P$ are more stringent under Assumption~\ref{as.str} alone than with also the use of Assumption~\ref{as.smooth}. This eliminates de-facto the requirement of knowing (or estimating) the bounds $C_0, C_1, C_2, C_3$, which may be difficult to do in practice. 


\subsection{Unconstrained algorithm convergence}

The results of Theorems~\ref{theo.constrained_0}-\ref{theo.constrained} can be tuned to the case of unconstrained problems as shown in the following. 

\vskip2mm

\begin{theorem}\label{theo.unconstrained_0}

{\bf
Consider the sequence $\{\x_k\}_{k \in \mathbb{N}}$ generated by the U-FOPC method, and let Assumption~\ref{as.str} hold true. Let $\x^*(t_k)$ be the optimizer of~\eqref{eq.problem} at time $t_k$ for $X = \mathbb{R}^n$. Define the following quantities 
\begin{equation}
\varrho_{\mathrm{P}} = \max\{|1-\alpha m|,|1- \alpha L|\}, \, \varrho_{\mathrm{C}} = \max\{|1-\beta m|,|1- \beta L|\},
\end{equation}
and choose stepsizes $\alpha$ and $\beta$ as
\begin{equation}
\alpha < 2/L, \quad \beta < 2/L.
\end{equation}
Further, select the number of correction steps $C$ in a way that  
\begin{equation}\label{eq.cvrt_u}
\tau_0 := \varrho_{\mathrm{C}}^{C} \left[\varrho_{\mathrm{P}}^{P} +  (\varrho_{\mathrm{P}}^{P}+1)\Big(1-\gamma +\gamma \frac{2L}{m}\Big)\right]  <1.
\end{equation} 
Then, the sequence $\{\|\x_k-\x^*(t_k)\|\}_{k \in \mathbb{N}}$ converges linearly with rate $\tau_0$ to an asymptotical error bound, and 
\begin{equation}\label{main.result2_0}
\limsup_{k\to\infty}\|\x_{k} - \x^*(t_k)\| = O(\varrho_{\mathrm{C}}^{C} h).
\end{equation}
}
\end{theorem}

\vskip2mm
\begin{IEEEproof}
See Appendix~\ref{ap.newton}, where we also derive all the constants in the right-hand side of~\eqref{main.result2_0}. The proofs leverage generalized equation theory and the special unconstrained nature of the problem.
\end{IEEEproof}

\vskip2mm

Similarly to Theorem~\ref{theo.constrained_0}, Theorem~\ref{theo.constrained_0} is a global $O(h)$ convergence result. From condition~\eqref{eq.cvrt_u}, it would seem that the best choice of $\gamma$ would be $\gamma = 0$; however this may not be the case when Assumption~\ref{as.smooth} holds. 

\vskip2mm

\begin{theorem}\label{theo.unconstrained}
{\bf
Consider the sequence $\{\x_k\}_{k \in \mathbb{N}}$ generated by the U-FOPC method. Let Assumptions~\ref{as.str}-\ref{as.smooth} hold true. Let $\x^*(t_k)$ be the optimizer of~\eqref{eq.problem} at time $t_k$ for $X = \mathbb{R}^n$. Define, 
\begin{equation}
\varrho_{\mathrm{P}} = \max\{|1-\alpha m|,|1- \alpha L|\}, \, \varrho_{\mathrm{C}} = \max\{|1-\beta m|,|1- \beta L|\},
\end{equation}
and set the stepsizes $\alpha$ and $\beta$ as
\begin{equation}
\alpha < 2/L, \quad \beta < 2/L.
\end{equation}
Let $\tau\in (0,1)$, $P$, $C$ be such that 
\begin{equation}
(1-\gamma)\varrho_{\mathrm{C}}^{C}(1 + \varrho_{\mathrm{P}}^{P}) + \varrho_{\mathrm{P}}^{P}\varrho_{\mathrm{C}}^{C} < \tau.
\end{equation}

There exist an upper bound on the sampling period $\bar{h}$ and a convergence region $\bar{R}$, such that if the sampling period is chosen as $h \leq \bar{h}$ and the initial optimality gap satisfy $\|\x_{0} - \x^*(t_0)\| \leq \bar{R}$, then the sequence $\{\|\x_k-\x^*(t_k)\|\}_{k \in \mathbb{N}}$ converges linearly with rate $\tau$ to an asymptotical error bound, and 
\begin{equation}\label{main.result2}
\limsup_{k\to\infty}\|\x_{k} - \x^*(t_k)\| = O(h^2\,\varrho_{\mathrm{C}}^{C}) + O(h\,\varrho_{\mathrm{P}}^{P}\varrho_{\mathrm{C}}^{C}).
\end{equation}
In addition, the bounds $\bar{h}$ and $\bar{R}$ are given as
\begin{align}\label{h.u}
\bar{h} &= \Big(\frac{\tau - \varrho_{\mathrm{C}}^{C}\varrho_{\mathrm{P}}^{P}}{\varrho_{\mathrm{C}}^{C}(\varrho_{\mathrm{P}}^{P} + 1)} - 1 + \gamma\Big)\Big(\frac{C_1 C_0}{m^2} + \frac{C_2}{m} \Big)^{-1}, \\ \bar{R} &= \frac{2 \,m}{\gamma C_1} \Big(\frac{C_1 C_0}{m^2} + \frac{C_2}{m}\Big)\, (\bar{h}-h).
\end{align}
}
\end{theorem}

\vskip2mm
\begin{IEEEproof}
See Appendix~\ref{ap.newton}, where we also derive all the constants in the right-hand side of~\eqref{main.result2}. The proofs leverage generalized equation theory and the special unconstrained nature of the problem. 
\end{IEEEproof}

\vskip2mm
Theorem~\ref{theo.unconstrained} is a local $O(h^2)$ convergence result. It can be seen as a generalization of Theorem~1 in~\cite{Paper1} and a special, yet more detailed, version of Theorem~\ref{theo.constrained}. In particular, by properly selecting the parameter $\gamma$, one obtains the results of both Theorem~\ref{theo.constrained} and~\cite{Paper1}; when  $\gamma = 1$, we obtain the  results of Theorem~\ref{theo.constrained} and this indicates that considering constrained problems does not add extra errors to the asymptotical bounds. On the other hand, if $\gamma = 0$, we are implicitly assuming that the prediction stage leads one to navigate into the iso-residual manifold;  since $\bar{R} \to \infty$, we also obtain global convergence. Furthermore, one can use $\gamma$ as a tuning mechanism to enlarge the convergence region. 

Finally, for $\gamma = 0$ and $P \to \infty$, then the results of theorem~\eqref{main.result2} boil down to those in~\cite{Paper1}. In particular, from~\eqref{h.u}, it can be seen  that the sampling period must satisfy the following relationship
\begin{equation}
\varrho_{\mathrm{C}}^{C}\Big[ 1 + h\, \Big(\frac{C_1 C_0}{m^2} + \frac{C_2}{m} \Big)\Big] < \tau < 1,
\end{equation}
which is the same requirement of~\cite[Theorem~1]{Paper1}.

\section{Numerical Experiments}\label{sec:num}

In this section, we report three numerical examples to showcase the performance of the proposed algorithms. First, we analyze the unconstrained case (Algorithm~\ref{algo_fou}), then we switch our attention to the constrained one (Algorithm~\ref{algo_fo}). Finally, we present a realistic simulation stemming from real-time control of energy resources. 

\subsection{Unconstrained example}


We use the same scalar example of~\cite[Section~IV.A]{Paper1}, and in particular we consider a time-varying cost function of the form
\begin{equation}\label{eq:scalar}
\min_{x \in \reals} f(x; t): = \frac{1}{2}\left(x - \cos(\omega t)\right)^2 + \kappa \log[1 + \exp(\mu x)]. 
\end{equation} 
The function in \eqref{eq:scalar} represents, for instance, the goal of staying close to a periodically varying trajectory plus a logistic term that penalizes large values of $x$. The terms $\omega$, $\kappa$, and $\mu$ are arbitrary nonnegative scalar parameters. In our experiments these parameters are set to $\omega = \pi/2$, $\kappa = 2$, and $\mu = 1.75$. The function $f(x; t)$  satisfies all the conditions in Assumptions \ref{as.str} and \ref{as.smooth}. In particular, $m = 1$ and $L = 2.53$.

\begin{figure}[t]
\centering
{
\psfrag{x}[c]{\footnotesize Sampling time $t_k$}
\psfrag{y}[c]{\footnotesize Tracking error $\|x_k - x^*(t_k)\|$}
\psfrag{runninggradientgradientgradientgradient}{\footnotesize \footnotesize U-FOPC, $P=0$, $C=3$, $\gamma = 0$}
\psfrag{gradient}{\footnotesize U-FOPC, $P=1$, $C=3$, $\gamma = 0$}
\psfrag{approxgradient}{\footnotesize U-FOPC, $P=1$, $C=3$, $\gamma = 1$}
\psfrag{newton3}{\footnotesize \footnotesize U-FOPC, $P=\infty$, $C=3$, $\gamma = 0$}
\psfrag{newton4}{\footnotesize \footnotesize U-FOPC, $P=\infty$, $C=3$, $\gamma = 1$}
\psfrag{newton}{\footnotesize\footnotesize U-FOPC, $P=3$, $C=3$, $\gamma = 0$}
\psfrag{newton2}{\footnotesize \footnotesize U-FOPC, $P=3$, $C=3$, $\gamma = 1$}
}
\includegraphics[width=0.5\textwidth]{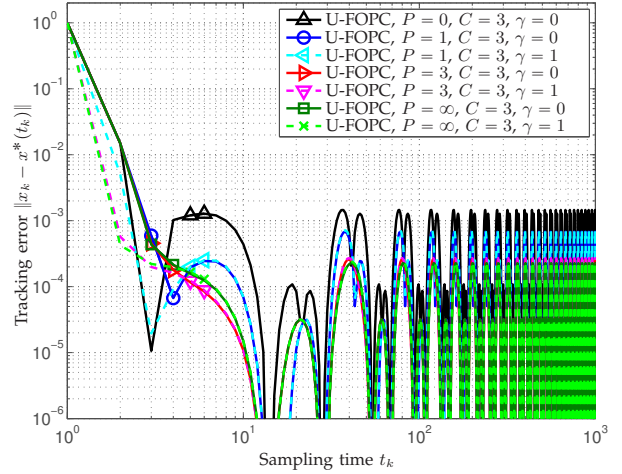}
\caption{Error with respect to the sampling time $t_k$ for different algorithms applied to the scalar problem~\eqref{eq:scalar}, with $h = 0.1$.}
\label{fig.centr1}
\end{figure}
\begin{figure}[t]
\centering
{
\psfrag{x}[c]{\footnotesize Sampling period $h$}
\psfrag{y}[c]{\footnotesize Asymptotical worst-case error}
\psfrag{runninggradientgradientgradientgradient}{\footnotesize \footnotesize U-FOPC, $P=0$, $C=3$, $\gamma = 0$}
\psfrag{gradient}{\footnotesize U-FOPC, $P=1$, $C=3$, $\gamma = 0$}
\psfrag{approxgradient}{\footnotesize U-FOPC, $P=1$, $C=3$, $\gamma = 1$}
\psfrag{newton3}{\footnotesize \footnotesize U-FOPC, $P=\infty$, $C=3$, $\gamma = 0$}
\psfrag{newton4}{\footnotesize \footnotesize U-FOPC, $P=\infty$, $C=3$, $\gamma = 1$}
\psfrag{newton}{\footnotesize U-FOPC, $P=3$, $C=3$, $\gamma = 0$}
\psfrag{newton2}{\footnotesize  U-FOPC, $P=3$, $C=3$, $\gamma = 1$}
\psfrag{oh}{\footnotesize $O(h)$}
\psfrag{oh2}{\footnotesize  $O(h^2)$}
}
\includegraphics[width=0.5\textwidth]{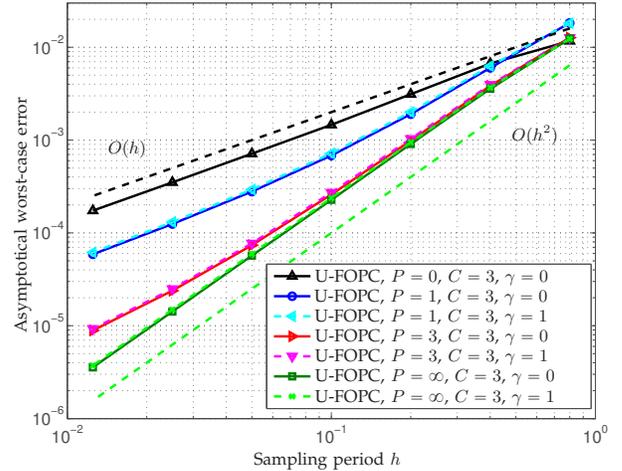}
\caption{Asymptotical worst case error floor with respect to the sampling time interval $h$ for different algorithms applied to the scalar problem~\eqref{eq:scalar}.}
\label{fig.centr2}
\end{figure}

We choose the constant stepsizes as $\alpha = \beta = 0.56< 2/L$ in the gradient methods stated in Algorithm~\ref{algo_fou} and initialize $x_0=0$ for all the algorithms. According to Theorem~\ref{theo.unconstrained_0}, for $\gamma = 1$, we need $C>2$ to ensure $\tau_0<1$ and therefore $O(h)$ convergence. In the case of Theorem~\ref{theo.unconstrained}, for $\gamma = 0$, $P = 1$ and $C=3$: we can set $\tau > 0.16$, the sampling period needs to be chosen as $h \leq \bar{h} =   4.2$, and the convergence region is the whole $\reals$, while for $\gamma =1$, $P = 1$ and $C=3$: we can set $\tau > 0.038$, $\bar{h} = 4.84$, and $\bar{R} = 12.7$ (notice that for greater $P$ or $C$ the requirements are less stringent).
 
In Figure~\ref{fig.centr1}, we plot the error $\|x_k - x^*(t_k)\|$ versus the discrete time $t_k$ for a sampling period of $h = 0.1$, for different schemes (we fix the correction steps at $C=3$ for all the algorithms and we study the behavior for different prediction steps $P$ and $\gamma$ values). Observe that the running gradient method~\cite{Popkov2005}, that is U-FOPC with no prediction and one correction step performs the worst, while U-FOPC with exact prediction (that is $P = \infty$) performs the best\footnote{The running gradient is a standard gradient descent with constant stepsize, where at each iteration one performs $C$ correction steps and the cost changes at each sampling time. }. Furthermore, we can notice that increasing $P$, one obtains better and better asymptotical error, with the drawback of an increased computational burden. Lastly, we notice how $\gamma =1$, that is having a Newton-like prediction step helps in achieving a faster convergence, and yet it appears that the asymptotical error is slightly greater than using $\gamma = 0$ (i.e., tangential prediction).  

The differences in performance can be also appreciated by varying $h$ and observing the worst case error floor size which is defined as 
$ 
\max_{k>\bar{k}} \{\|x_k - x^*(t_k)\|\},
$
where $\bar{k} = 10^4$ in the simulations. Figure~\ref{fig.centr2} illustrates the error as a function of $h$. The performance differences between the proposed methods that may be observed here corroborate the differences seen in Figure \ref{fig.centr1}. In particular, the U-FOPC method with $P=0, C=1$ achieves the largest worst case error bound, while exact prediction with $\gamma = 0$ attains the best case error bound. Notice also the dashed lines displaying the theoretical performance of $O(h)$, $O(h^2)$, and the fact that for any finite $P$ there is a trade-off between computational complexity and asymptotical error\footnote{The simulation code is made available online to further analyze how different parameters can affect the algorithm trade-offs.}. 

\subsection{Constrained example}\label{sec:numconstr}

We consider here a mid-size optimization problem consisting of $n = 1000$ scalar variables. The cost function we study has the form
\begin{equation}
f(\x; t) = \frac{1}{2}\|\x + {\bf 1}_n\|^2_{\Q} + \sum_{i=1}^n \kappa_i \sin^2(\omega t + \varphi_i)\,\exp(\mu (x_{(i)}-2)^2),
\end{equation}
where we have defined ${\bf 1}_n$ as the column vector of all ones of dimension $n$, while $x_{(i)}$ is the $i$-th component of $\x\in \reals^n$. In addition, the matrix $\Q$ is chosen as $\Q = {\bf I}_n + \vv\vv^\transp/n$ with $\vv$ being a vector randomly generated by a normal distribution of mean $0$ and variance $1$, $\kappa_i \sim \mathcal{U}_{[0,1]} $, $\omega = 0.1\pi $, $\varphi \sim \mathcal{N}(0,\pi)$, and $\mu = 0.25$. 

We study the time-varying problem
\begin{equation}\label{prob.vector}
\minimize_{\x \in [0,0.4]^n} \, f(\x; t).
\end{equation}
We notice that the cost function $f$ verifies the Assumptions~\ref{as.str}-\ref{as.smooth} on on $[0,0.4]^n$, which is our optimization set (even though it would not satisfy them on the whole $\reals^n$). In particular, $m = 1$ and $L = 6.07$. 

One could run a similar analysis as the unconstrained example, however here we focus on realistic run-time constraints. Every time a new function is available, a number of correction steps are performed. The number depends on how fast we need the corrected variable to be available and the computational time necessary to compute the gradient and perform the correction step. We fix at $r_1 h$, with $r_1<1$ the time allocated for the correction steps, while $t_{\textrm{C}}$ is the time to perform one correction step. For the above considerations, we can afford to run
\begin{equation}\label{opt.c}
C = \lfloor r_1 h/ t_{\textrm{C}}\rfloor,
\end{equation}  
correction steps. After the corrected variable is available, one can use it for the decision making process (which may require extra time to be performed). For the time-varying algorithm perspective, one can use the variable to either run $P$ gradient prediction, or $C'$ extra correction steps (to improve the corrected variable for having a better starting point when a new function becomes available). Fix at $r_2 h$, with $r_2<1$ the time allocated for the prediction (or extra correction) steps. The affordable number of prediction steps can be determined considering that $P$ prediction steps require a time equal to $\bar{t} + P t_{\textrm{P}}$, where $\bar{t}$ is the time required to evaluate the Hessian, gradient, and time derivative of the gradient, while $t_{\textrm{P}}$ is the time to perform one prediction calculation. Thus, 
\begin{equation}\label{opt.p}
P = \lfloor (\,r_2 h - \bar{t}\,)/ t_{\textrm{P}}\rfloor.
\end{equation}  
The affordable extra correction steps $C'$ can be computed as in~\eqref{opt.c}, substituting $r_1$ with $r_2$. 

In the simulation example, we choose $r_1 = r_2 = 0.5$, while by running the experiments on a $1.8$~GHz Intel Core i5, we empirically fix $t_{\textrm{C}} = .76$~ms, $\bar{t} = 10$~ms, $t_{\textrm{P}} = .62$~ms. Note that the time that would be needed to solve the prediction step exactly (by solving a quadratic program) is $190$~ms, which is not affordable in the considered sampling period range.

In addition, we consider the situation in which one can use the whole sampling period to do correction, that is $r_1 = 1$, while $r_2 = 0$, and we call this case \emph{total correction}. This situation is particularly interesting when one has to make a choice whether to stop the correction steps to perform prediction, or to continue to do correction steps till a new function evaluation becomes available. Note that the correction+extra correction strategy is different from the total correction one, since the error is computed with the corrected variable (which is used for the decision making process), that is after $r_1 h$.

\begin{figure}[t]
\centering
{
\psfrag{x}[c]{\footnotesize Sampling period $h$ in [ms]}
\psfrag{y}[c]{\footnotesize Asymptotical worst-case error}
\psfrag{Ecorr}{\footnotesize Correction + Extra correction, $r_1$ = $0.5$, $r_2$ = $0.5$}
\psfrag{CorrTotCorrTotCorrTotCorrTot}{\footnotesize Total correction, $r_1$ = $1$, $r_2$ = $0$}
\psfrag{Pred}{\footnotesize C-FOPC Algorithm, $r_1$ = $0.5$, $r_2$ = $0.5$}
}
\includegraphics[width=0.5\textwidth]{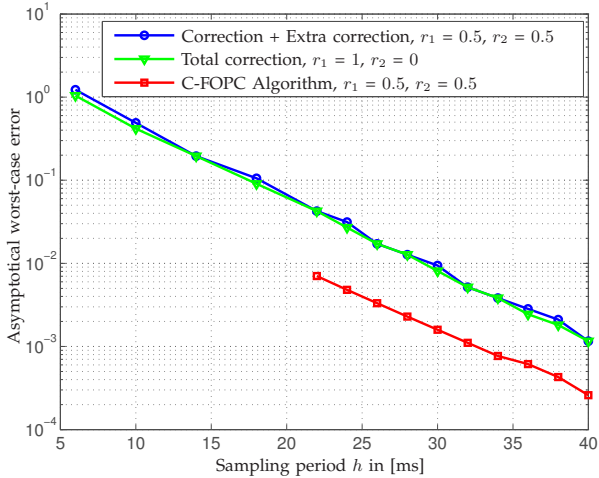}
\caption{Asymptotical worst case error floor with respect to the sampling time interval $h$ for different algorithms applied to~\eqref{prob.vector}}
\label{fig.centr3}
\end{figure}

In Figure~\ref{fig.centr3}, we report the asymptotical worst-case error w.r.t the sampling period for the three considered cases (correction+extra correction, total correction, and prediction-correction, i.e. the C-FOPC algorithm)\footnote{Note that both the correction+extra correction and the total correction strategies are in fact (running) projected gradient descent algorithms.}, while the number of prediction steps and correction steps are optimized via the available resources as in Eq.s~\eqref{opt.c}-\eqref{opt.p}. With the simulation parameters, for $h = 6$~ms, we can perform $C = C' = 3$ steps of correction and extra correction, or $C = 7$ steps of (total) correction. For $h = 40$~ms, these values are $C = C' = 26$ and $C = 52$, respectively. For the prediction-correction strategy, for $h = 22$~ms, then $C = 14$ and $P = 1$, while for $h = 40$~ms, $C = 26$ and $P = 16$. 

For sampling times below $22$~ms, prediction cannot be performed due to time constraints. For sampling periods greater or equal than $22$~ms, prediction can be performed and for $h = [22, 40]$~ms, then $\bar{h} = [90, 370]$~ms, and $\bar{R} = [.13, .68]$. We see clearly that, in this simulation example, if prediction is affordable, the prediction-correction strategy, that is our C-FOPC algorithm is to be preferred to traditional correction-only schemes, since it achieves a lower asymptotical worst-case error. We notice that this error is lower by an half order of magnitude, while the error of the correction-extra correction and the total correction strategy are practically the same. For completeness, we report that $\x_0$ is chosen to be zero.


The result is quite remarkable, telling that performing Newton-like prediction steps on a fixed (Hessian, gradient, time derivative) triple can be computationally much more interesting that performing correction steps on a varying (i.e., re-updated) gradient.  

In Figure~\ref{fig.centr4}, we report the time trajectories of a number of variables for the three strategies to appreciate how the constraints are in fact active. 

\begin{figure}[t]
\centering
{
\psfrag{x}[c]{\footnotesize Time [s]}
\psfrag{y}[c]{\footnotesize Trajectories $\{\x_k\}$ and $\x^*(t)$}
\psfrag{xopt}{\footnotesize Optimal trajectory $\x^*(t)$}
\psfrag{Ecorr}{\footnotesize Correction + Extra correction, $r_1$ = $0.5$, $r_2$ = $0.5$}
\psfrag{CorrTotCorrTotCorrTotCorrTot}{\footnotesize Total correction, $r_1$ = $1$, $r_2$ = $0$}
\psfrag{Pred}{\footnotesize C-FOPC Algorithm, $r_1$ = $0.5$, $r_2$ = $0.5$}
}
\includegraphics[width=0.5\textwidth]{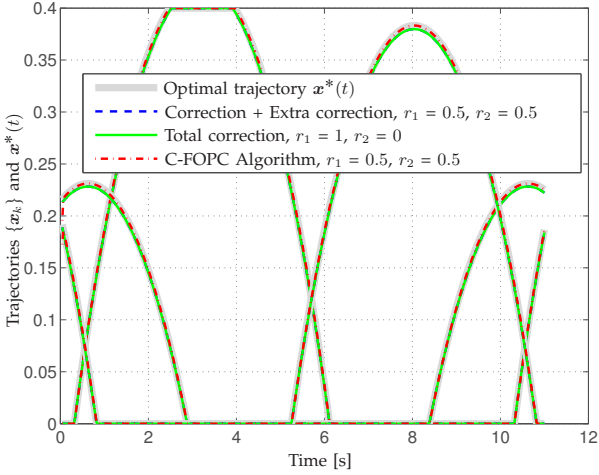}
\caption{Trajectories of $\{\x_k\}$ and optimal trajectory $\x^*(t)$ for different algorithms applied to~\eqref{prob.vector}, with $h = 22$~ms.}
\label{fig.centr4}
\end{figure}

\subsection{An example of realistic application}

As an example of realistic application, we consider the problem of optimizing  the operation of aggregations  of  energy resources connected to a power distribution system -- a research task that has gained significant interest from academic and industrial sectors because of the increased  deployment of distributed energy resources (DERs).  The proposed algorithm can be implemented to enable a \emph{real-time} optimization of the DER operation, where the term ``real-time'' refers to a setting where the power setpoints of the DERs are updated on a second or subsecond time scale based on streaming measurements, to maximize operational objectives while coping with variability of ambient conditions and non-controllable energy assets~\cite{DallAnese2016, LowOnlineOPF, AndreayOnlineOpt}.

Consider then a distribution network consisting of one slack bus and $N$ nodes with controllable DERs. Let $s_n = p_n + j q_n \in \mathbb{C}$ denote the net injected complex power at node $j$,  and let $p_0  \in  \mathbb{R}$ denote  the active power flow at a point of connection of the distribution system with the rest of the grid.  For future developments, let $\x_{n} := [\Re\{s_n\}, \Im\{s_n\}]^\transp$ collect the active and reactive setpoints of the DER at node $n$. To further simplify the notation, define the vector $\x := [\x_{1}^{\transp}, \ldots, \x_N^{\transp}]^{\transp}$. To facilitate the development of computationally-affordable algorithms, we postulate the following approximate linear relationship between $\x$ and the active powers $p_0$:  
\begin{align} 
\tilde{p}_0(\x, t) & =  \m^{\transp}(t) [\x - \x_\ell(t)] + n(t) \nonumber \\
& := \sum_{n = 1}^N \m^{\transp}_n(t)[\x_{n} - \x_{n, \ell}(t) ] + \eta(t) \label{eqn:lin_s0}
\end{align}
where $\x_{n,\ell}(t)$ denotes the non-controllable loads and the model parameters $\{\bm_n(t) \in \mathbb{R}^2\}$ and $\eta(t)  \in \mathbb{R}$ can be obtained as shown in, e.g.,~\cite{DallAnese2016,bolognani2015linear}. These model parameters can be either fixed or time-varying.

Consider a discrete-time operational setting where the setpoints of the DERs are updated at time instants $t_k = h k$, with $k \in \mathbb{N}$ and $h > 0$ based on the specific implementation requirements (e.g., second, sub-second,  or a few seconds), and formulate the following optimization problem per time $k$:
\begin{subequations} 
\label{eq:powerproblem}
\begin{align} 
\min_{\substack{\x}\in\mathbb{R}^{2n}} & \, \sum_{n=1}^N  C_{n}(\x_{n}; t_k) + \frac{\gamma}{2} \left (p_0^{\mathrm{set}}(t_k) - \tilde{p}_0(\x, t_k)\right )^2 \\
\textrm{s.~t.:~} & \quad \x_{n} \in Y_{n}, \, \forall \,\,\, n = 1, \dots, N \label{eqn:constr_Y_exact} 
 \end{align}
\end{subequations}
where  $p_0^\mathrm{set}(t_k) $ is a setpoint for  the aggregate power at the point of connection~\cite{AndreayOnlineOpt},  $\gamma > 0$ is a design parameter that influences the ability to track the reference signal $p_0^\mathrm{set}(t_k) $, $C_n(\x_{n}; t_k):\mathbb{R}^2\times\mathbb{R}_{+} \to \mathbb{R}$ is a local cost function, and $Y_{n}$ specifies the operating region for DER $n$. For example, $Y_{n}$ can be given by
\begin{align} 
Y_{n} = \{p_n, q_n: \underline{p}_n \leq p_n \leq \overline{p}_n, \underline{q}_n \leq q_n \leq \overline{q}_n\}
 \end{align}
with $ \underline{p}_n$, $\overline{p}_n$, $ \underline{q}_n$, and $\overline{q}_n$ given limits. This is the case, for example, for variable-speed drives, natural-gas generators, electric vehicle chargers, and renewable sources of energy operating with (fixed) headroom from the maximum available power.  Problem~\eqref{eq:powerproblem} is pertinent for settings where a number of customers aim at maximizing given performance objectives related to the DERs (e.g., minimization of net electricity payment), while partaking in grid services;  specifically, customers may be incentivized by utility or aggregators to adjust the DERs' output powers so that the aggregate power $\tilde{p}_0(\x; t_k)$ can be driven to the setpoint $p_0^\mathrm{set}(t_k)$.

For the numerical experiments, we utilize real data belonging to residential loads from a neighborhood located in the Sacramento metro area~\cite{Bank13}. The data have a granularity of 1 second and they are utilized to populate the IEEE 37-node test feeder (see e.g.,~\cite{DallAnese2016}). For simplicity, we consider a single-phase equivalent of the feeder. The model parameters $\{\m_n(t) \in \mathbb{R}^2\}$ and $\eta(t)  \in \mathbb{R}$ are computed  as shown in~\cite{DallAnese2016}.  With this setting, the non-controllable portion of $\tilde{p}_0(\x; t_k)$ in~\eqref{eq:powerproblem}, which defined as 
\begin{align} 
a(t_k) := - \sum_{n = 1}^N \m^{\transp}_n(t_k)\x_{n, \ell}(t_k) + \eta(t_k) 
\end{align}
varies over time as illustrated in Figure~\ref{fig:mov_a}. The trajectory for the setpoint $p_0^\mathrm{set}(t_k)$ is also provided in Fig.~\ref{fig:mov_a}; in the considered setting, the specified setpoints promote a smoother variability of the aggregate load $\tilde{p}_0(\x; t)$ over time. 

With reference to the feeder in~\cite{DallAnese2016}, we assume that controllable DERs are located at nodes 4, 7, 10, 13, 17, 20, 22, 26, 30, and 35,  and it is assumed that the real power $p_n$ can be controlled within the interval $ \underline{p}_n = -50$ kW and $\overline{p}_n = 50$ kW. The cost functions $C_n(\x_n;t_k) = p_n^2$ and $\gamma = 2$. This setting is representative of the case where DERs can deviate their active powers from a nominal operating point to provide services to the grid, but their adjustment is limited within $[\underline{p}_n, \overline{p}_n]$. In addition, $\m_n(t_k)$ is considered constant in time (as often varies at slower time-scales than $a(t_k)$, typically minutes).

To track the solution trajectory $\x^*(t_k)$ of~\eqref{eq:powerproblem} while $k$ increases, we run both a running projected gradient algorithm, i.e., Algorithm~\ref{algo_fo} with $P=0$, and our prediction-correction Algorithm~\ref{algo_fo}.  The prediction step is computed by approximating the time-derivative of the gradient as explained in~\cite{Paper1}, which does not affect convergence or the rates. 

For the simulation runs, the functional parameters are $m = 1$ and $L = 21$, while since the cost function is a time-varying quadratic function with constant Hessian, $C_1 = 0$, $C_2=0$, the convergence attraction region is the whole space and there is no upper bound on the sampling period (see Theorem~\ref{theo.constrained}). 
We choose the stepsizes $\alpha$ and $\beta$ as $\alpha = \beta = 0.0048$, such that $\varrho_{\textrm{P}} = \varrho_{\textrm{C}}<1 $, while $h = 1$~s, and $\x_0 = \bf{0}$.

We keep the computational time fixed in our comparison of running projected gradient and prediction-correction algorithm, so that we could identify better the contribution of the prediction step. With the same notation and nomenclature of Section~\ref{sec:numconstr}, we look at a total-correction strategy for the running projected gradient with $C = 3$. Each correction step takes $76$~$\mu$s on a $1.8$~GHz Intel Core i5 (and $3$ steps, $228$~$\mu$s). For the prediction-correction algorithm, we fix $C = 1$, while we choose $P = 2$, such that the total computational time is $206$~$\mu$s, which is similar to the running strategy. 

We report the error $\|\x_{k} - \x^*(t_k)\|$ of the two strategies in time in Figure~\ref{fig:fig6} (for readability purposes only for the time frame $12$:$38$-$12$:$48$, but in fact it is qualitatively the same for all the other time frames), and in particular the averaged asymptotical error (averaged in the time frame $12$:$10$-$13$:$00$) of the running method ($P=0$) is $3.73$~kW, while the one of the prediction-correction algorithm is $3.14$~kW. As captured in Figure~\ref{fig:fig6}, the prediction-correction algorithm seems to be more reactive to changes, responding more quickly, which is expected. This further corroborates the usefulness of the prediction step, hence of the prediction-correction methodology, even when real noisy time-varying data is used, and the time derivative of the gradient can only be approximated. 

\begin{figure}[t]
\centering
{
\psfrag{x}[c]{\footnotesize Time}
\psfrag{y}[c]{\footnotesize Aggregate active power [kW]}
\psfrag{Non-controllable aggregate load}{\footnotesize Non-controllable aggr. load}
\psfrag{Setpoint}{\footnotesize Setpoint}
}
\includegraphics[width=0.5\textwidth]{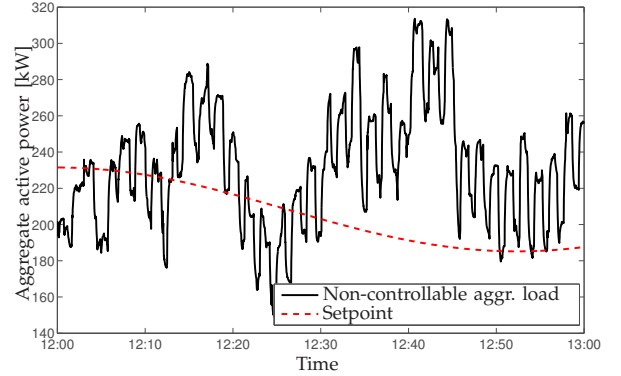}
\caption{Aggregate effect of the non-controllable residential load on the power at the point of common coupling $\tilde{p}_0(\x; t_k)$ [kW] in the considered simulation setting. The red trajectory represents the prescribed setpoint $p_0^\mathrm{set}(t_k)$. }
\label{fig:mov_a}
\end{figure}

\begin{figure}[t]
\centering
{
\psfrag{x}[c]{\footnotesize Time}
\psfrag{y}[c]{\footnotesize Tracking error $\|x_k - x^*(t_k)\|$}
\psfrag{aaaaaaaaaaaaaaaaaaaaaaaaaaaaaaaa}{\footnotesize \footnotesize C-FOPC, $P=0$, $C=3$}
\psfrag{b}{\footnotesize \footnotesize C-FOPC, $P=2$, $C=1$}
}
\includegraphics[width=0.5\textwidth]{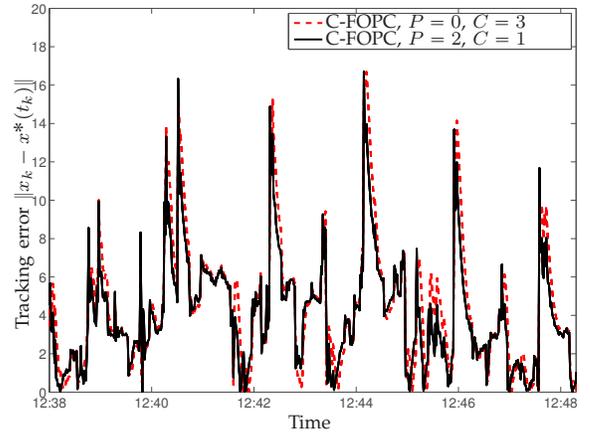}
\caption{Tracking error in [kW] for the two considered strategies in the time frame $12$:$38$-$12$:$48$. }
\label{fig:fig6}
\end{figure}

\section{Conclusions}\label{sec:concl}

We have proposed first-order algorithms to find and track the solution trajectory of strongly convex, strongly smooth constrained time-varying optimization problems. These algorithms rely on a discrete-time prediction-correction strategy, by which at each sampling time, the decision variables are corrected through one or multiple projected gradient steps, and then used to predict the next decision variables via successive projected gradient steps on a suitably defined quadratic program. The proposed algorithms exhibit better asymptotical accuracy than state-of-the-art correction-only schemes, even when computational complexity issues are taken into account. 


%

\appendices



\section{Preliminaries}

We begin the convergence analysis by deriving an upper bound on the norm of the approximation error $\mathbold{\Delta}_k\in \reals^n$ incurred by the Taylor expansion in~\eqref{ger0}. In particular, given  the optimal solution $\x^*(t_k)$ at $t_k$, compute the optimal prediction step via the Taylor approximation~\eqref{ger0} and indicate the optimal prediction as $\x^*_{k+1|k}$. The objective is to bound the error:
\begin{align}\label{eq.def.delta}
\mathbold{\Delta}_k : = \x^*_{k+1|k} - \x^*(t_{k+1}),
\end{align}
which is committed when $\x^*(t_{k+1})$ is replaced by $\x^*_{k+1|k}$ (here we use the superscript $*$ to indicate that we compute the prediction from $\x^*(t_k)$ and not any $\x_k$).


\vspace{1mm}
\begin{proposition}\label{prop.err}
{\bf
The error norm $\|\mathbold{\Delta}_k\|$ is upper bounded as follows. 

\emph{Case a)} Under the sole Assumption~\ref{as.str}
\begin{equation}\label{prop_claim_err_bound_0}
\|\mathbold{\Delta}_k\| \leq 2 h \frac{C_0}{m} (1+L/m) =: \Delta_1 = O(h).
\end{equation}

\emph{Case b)} Under both Assumption~\ref{as.str} and~\ref{as.smooth}
\begin{equation}\label{prop_claim_err_bound}
\|\mathbold{\Delta}_k\| \leq \frac{h^2}{2} \left[\frac{C_0^2 C_1}{m^3} + \frac{2 C_0 C_2}{m^2} + \frac{C_3}{m}\right] =: \Delta_2 = O(h^2).
\end{equation}
}
\end{proposition}

\begin{IEEEproof}
Let us start by simplifying the notation. Define
\begin{align}\label{simplif}
\nabla_{\x} f_i = \nabla_{\x} f (\x^*(t_{k+i}); t_{k+i}), \quad & \Q_i = \nabla_{\x\x} f (\x^*(t_{k+i}); t_{k+i}) \\ \label{simplif22}
\c_i = \nabla_{t\x} f (\x^*(t_{k+i}); t_{k+i}),\quad & \x_i = \x^*(t_{k+i}), \, \x = \x^*_{k+1|k}.
\end{align}
With this notation in place, $\mathbold{\Delta}_k = \x - \x_{1}$ [Cf.~\eqref{eq.def.delta}]. In addition, $\x$ is computed by the generalized equation~\eqref{ger0}, 
\begin{equation}\label{proofA.dummy1}
\nabla_{\x} f_0 + \Q_0 (\x - \x_0) + h\, \c_0 + N_X(\x) \ni \mathbf{0},
\end{equation}
while $\x_1$ is the solution of 
\begin{equation}\label{proofA.dummy2}
\nabla_{\x} f_1 + N_X(\x_1) \ni \mathbf{0}.
\end{equation}

Define the supporting functions, 
\begin{align}\label{proofA.dummy3}
G(\y) &= \nabla_{\x} f_1 + \Q_1(\y - \x_1) + N_X(\y)\\
g(\y) &= \nabla_{\x} f_0 + \Q_0 (\y- \x_0) + h\, \c_0 + \nonumber \\ & \hskip3cm - [\nabla_{\x}f_1 + \Q_1 (\y - \x_1)].  
\end{align}
These two functions allows one to rewrite~\eqref{proofA.dummy1} as 
\begin{equation}\label{proofA.dummy4}
(g + G)(\x) \ni \mathbf{0}.
\end{equation}
It is also true that $(g + G)(\x_1) \ni g(\x_1)$, since for optimality $G(\x_1) \ni \mathbf{0}$. Define $F(\y) = (g + G)(\y)$, and consider the parametric generalized equation $F(\y) + \p \ni \mathbf{0}$. Under the Assumptions~\ref{as.str}-\ref{as.smooth}, as for \cite[Theorem~2F.9]{Dontchev2009}, we have that the solution mapping $\p \mapsto \y(\p)$ for the generalized equation $F(\y) + \p \ni \mathbf{0}$ is every-where single valued and Lipschitz continuous with constant $m^{-1}$, that is 
\begin{equation}\label{proofA.dummy5}
\|\y(\p) - \y(\p')\| \leq \frac{1}{m} \|\p - \p'\|.
\end{equation}

We can set $\p = \mathbf{0}$ and $\p' = -g(\x_1)$, which leads to 
\begin{equation}\label{proofA.dummy6}
\|\x - \x_1\| = \|\mathbold{\Delta}_k\| \leq \frac{1}{m} \|g(\x_1)\|.
\end{equation}

We proceed now to bound $\|g(\x_1)\|$. We can write $g(\x_1)$ as
\begin{equation}\label{eq-dummy}
g(\x_1) = \nabla_{\x} f_0 + \Q_0 (\x_1- \x_0) + h\, \c_0 - \nabla_{\x} f_1.
\end{equation}


\subsection{Case a)}

Given the expression of $g(\x_1)$, we can bound its norm via Assumption~\ref{as.str}
\begin{align}\label{eq.bound1}
\|g(\x_1)\| &\leq \|\nabla_{\x} f_0 - \nabla_{\x} f_1\| + \|\Q_0\|\|\x_1- \x_0\| + h\, \|\c_0\| \nonumber \\
&\leq \|\nabla_{\x} f_0 - \nabla_{\x} f_1\| + L \|\x_1- \x_0\| + h C_0,
\end{align}
where we have used the upper bound on the Hessian of the cost function. Furthermore, 
\begin{multline}\label{eq.bound2}
\!\!\!\|\nabla_{\x} f_0 - \nabla_{\x} f_1\|  = \|\nabla_{\x} f(\x^*(t_k);t_k) - \nabla_{\x} f(\x^*(t_{k+1});t_{k+1})\| \\
\leq \|\nabla_{\x} f(\x^*(t_k);t_{k+1}) - \nabla_{\x} f(\x^*(t_{k+1});t_{k+1})\| +\\ \|\nabla_{\x} f(\x^*(t_k);t_{k}) - \nabla_{\x} f(\x^*(t_{k});t_{k+1}) \|
\end{multline}
We can bound the first term in the right hand side with the Lipschitz property of the gradient of the cost function as
\begin{equation}\label{eq.bound3}
\|\nabla_{\x} f(\x^*(t_k);t_{k+1}) - \nabla_{\x} f(\x^*(t_{k+1});t_{k+1})\| \leq L \|\x_1 - \x_0\|,
\end{equation}
while the second part can be bounded by the upper bound on $\nabla_{t\x}f$, as
\begin{multline}\label{eq.bound4}
\!\!\|\nabla_{\x} f(\x^*(t_k);t_{k}) - \nabla_{\x} f(\x^*(t_{k});t_{k+1}) \| \!= \! \|\nabla_{\x} f(\x^*(t_k);t_{k})  \\-\nabla_{\x} f(\x^*(t_{k});t_{k}) - h\, \nabla_{t\x}f(\x^*(t_{k});\tau)\| \leq h C_0,
\end{multline}
where we have used the mean value theorem and $\tau \in [t_k, t_{k+1}]$. By putting Eq.s~\eqref{eq.bound1}, \eqref{eq.bound2}, \eqref{eq.bound3}, and~\eqref{eq.bound4} together, one arrives at
\begin{equation}\label{eq-dummy_1}
\|g(\x_1)\| \leq 2 L \|\x_1- \x_0\| + 2 h\, C_0,
\end{equation}
and by using~\eqref{proofA.dummy6} as well as the bound on $\|\x_1 - \x_0\|$ in~\eqref{eq.lip}, the claim~\eqref{prop_claim_err_bound_0} follows.


\subsection{Case b)}
By allowing higher order smoothness, the bound on $g(\x_1)$ can be improved as follows. Consider~\eqref{eq-dummy}: it is nothing else but the error of the truncated Taylor expansion of $\nabla_{\x} f_1$:
\begin{equation}
\nabla_{\x} f_1 - \nabla_{\x} f_0 =  \Q_0 (\x_1-\x_0) + h \c_0 +\mathbold{\epsilon},
\end{equation} 
where the error $\mathbold{\epsilon}$ can be bounded as
\begin{multline}
\|\mathbold{\epsilon}\| \leq \frac{1}{2}\Big(\|\nabla_{\x\x\x} f\| \|\x_1- \x_0\|^2 + h\, \|\nabla_{t\x\x} f\| \|\x_1- \x_0\| +\\  h\, \|\nabla_{\x t\x} f\| \|\x_1- \x_0\|  + {h^2} \|\nabla_{t t\x} f\|\Big),
\end{multline}
and by using the upper bounds in Assumption~\ref{as.smooth}, 
\begin{equation}
\|g(\x_1)\| \leq \frac{1}{2}\, C_1 \|\x_1 - \x_0\|^2 + h\, C_2 \|\x_1 - \x_0\| + \frac{1}{2} h^2\, C_3. 
\end{equation}
By using the bound~\eqref{eq.lip} on the variability of the optimizers $\x_1$ and $\x_0$, then 
\begin{equation}
\|g(\x_1)\| \leq h^2\,C_1 \frac{C_0^2}{2 m^2} + h^2\, C_2 \frac{C_0}{m} + \frac{1}{2} h^2\, C_3, 
\end{equation}
and by combining this bound with~\eqref{proofA.dummy6}, the claim~\eqref{prop_claim_err_bound} follows. 
\end{IEEEproof}

\section{Proof of Theorems \ref{theo.constrained_0} and \ref{theo.constrained}}\label{ap.grad}

We divide the proof in different steps. Step 1: we bound the prediction error in Propositions~\ref{prop.intheo.1} and~\ref{prop.intheo.2}; Step 2: we bound the correction error; Step 3: we put the previous steps together and derive the convergence requirements and results. As done for Proposition~\ref{prop.err}, Case a) will refer to Theorem~\ref{theo.constrained_0}, while Case b) to Theorem~\ref{theo.constrained}.

\vskip2mm 

\noindent \textbf{Prediction error.} The optimal prediction error, i.e., the distance between the optimal predicted variable $\x_{k+1|k}$ and the optimizer at time step $t_{k+1}$, $\x^*(t_{k+1})$ can be bounded as the following proposition. 
\vskip2mm

\begin{proposition}\label{prop.intheo.1}
{\bf
The following facts hold true. 
\begin{enumerate}
\item[{Case a)}] Under the same assumptions and notation of Theorem~\ref{theo.constrained_0}, we have that
\begin{equation}\label{prop_inside_0}
\|\x_{k+1|k} - \x^*(t_{k+1})\| \leq \frac{2 L}{m}\|\x_k- \x^*(t_k)\|+ \Delta_1.
\end{equation}
\item[{Case b)}] Under the same assumptions and notation of Theorem~\ref{theo.constrained}, we have that
\begin{multline}\label{prop_inside}
\|\x_{k+1|k} - \x^*(t_{k+1})\| \leq \frac{C_1}{2\,m}\|\x_k- \x^*(t_k)\|^2 +\\ h\,\Big(\frac{C_1 C_0}{m^2} + \frac{C_2}{m}  \Big)\|\x_k- \x^*(t_k)\|+ \Delta_2.
\end{multline}
\end{enumerate}
}
\end{proposition}
\vskip2mm
\begin{IEEEproof}
We proceed as in the proof of Proposition~\ref{prop.err}. We use similar simplifications of~\eqref{simplif}, as
\begin{subequations}\label{simplif2}
\begin{align}
\nabla_{\x} f_k = \nabla_{\x} f (\x_k; t_{k}), \quad & \Q_k = \nabla_{\x\x} f \x_k; t_{k}) \\
\c_k = \nabla_{t\x} f (\x_k; t_{k}),\quad & \x = \x_{k+1|k}.
\end{align}
\end{subequations}
while $\nabla_{\x} f_1$ and $\x_1$ are defined just as in~\eqref{simplif}. The error $\|\x_{k+1|k} - \x^*(t_{k+1})\|$ is now $\x - \x_1$. 

The vector $\x$ is computed by the generalized equation~\eqref{ger0}, 
\begin{equation}\label{proofB.dummy1}
\nabla_{\x} f_k + \Q_k (\x - \x_k) + h\, \c_k + N_X(\x) \ni \mathbf{0},
\end{equation}
while $\x_1$ is the solution of~\eqref{proofA.dummy2}.

Define the supporting functions, 
\begin{align}\label{proofB.dummy3}
G(\y) &= \nabla_{\x} f_1 + \Q_1(\y - \x_1) + N_X(\y)\\
g(\y) &= \nabla_{\x} f_k + \Q_k (\y- \x_k) + h\, \c_k + \nonumber \\ & \hskip3cm - [\nabla_{\x}f_1 + \Q_1 (\y - \x_1)].  
\end{align}
These two functions allows one to rewrite~\eqref{proofB.dummy1} as 
\begin{equation}\label{proofB.dummy4}
(g + G)(\x) \ni \mathbf{0}.
\end{equation}
It is also true that $(g + G)(\x_1) \ni g(\x_1)$, since for optimality $G(\x_1) \ni \mathbf{0}$. Define $F(\y) = (g + G)(\y)$, and consider the parametric generalized equation $F(\y) + \p \ni \mathbf{0}$. Due to Assumptions~\ref{as.str}-\ref{as.smooth}, and due to \cite[Theorem~2F.9]{Dontchev2009}, we have that the solution mapping $\p \mapsto \y(\p)$ of the generalized equation $F(\y) + \p \ni \mathbf{0}$ is every-where single valued and Lipschitz continuous as
\begin{equation}\label{proofB.dummy5}
\|\y(\p) - \y(\p')\| \leq \frac{1}{m} \|\p - \p'\|.
\end{equation}

Set $\p = \mathbf{0}$ and $\p' = -g(\x_1)$, then, 
\begin{equation}\label{proofB.dummy6}
\|\x - \x_1\| \leq \frac{1}{m} \|g(\x_1)\|.
\end{equation}

We proceed now to bound $\|g(\x_1)\|$. We can write $g(\x_1)$ as
\begin{equation}\label{g1_ex}
g(\x_1) = \nabla_{\x} f_k + \Q_k (\x_1- \x_k) + h\, \c_k - \nabla_{\x} f_1,
\end{equation}
As in Proposition~\ref{prop.err}, we have two cases: 
\begin{align}
\textrm{Case a)} \quad \|g(\x_1)\| \leq  & 2 L\|\x_1 - \x_k\| + 2 h\, C_0, \label{g1_bound_0} \\ 
%
%
\textrm{Case b)} \quad  \|g(\x_1)\| \leq  &  \frac{1}{2}\, C_1 \|\x_1 - \x_k\|^2 + \nonumber \\ & \qquad h\, C_2 \|\x_1 - \x_k\| + \frac{1}{2} h^2\, C_3.  \label{g1_bound}
\end{align}
Since $\|\x_1 - \x_k\| \leq \|\x_1 - \x^*(t_k)\| + \|\x^*(t_k) - \x_k\|$, and we can bound the first term of the right-hand side by using~\eqref{eq.lip}, then 
\begin{align}
\textrm{Case a)}\, \|g(\x_1)\| \leq & 2 L\|\x^*(t_k) - \x_k\|\! + \!2 h C_0 (1 \!+\! {L}/m), \\
\textrm{Case b)}\, \|g(\x_1)\| \leq & h^2\,C_1 \frac{C_0^2}{2 m^2} + h^2\, C_2 \frac{C_0}{2 m} +  \frac{1}{2} h^2\, C_3 + \nonumber \\ & 
 \hspace*{-2.5cm}\frac{1}{2} C_1 \|\x^*(t_k) \!-\! \x_k\|^2 \!+\! h \Big(\frac{C_1 C_0}{m} \!+\! C_2 \Big) \|\x^*(t_k) \!-\! \x_k\|.
\end{align}
By combining these two last bounds with~\eqref{proofB.dummy6}, the claims~\eqref{prop_inside_0}-\eqref{prop_inside} follow.
\end{IEEEproof}

\vskip2mm

On the other hand, the distance between the approximate prediction $\tilde{\x}_{k+1|k}$ and the optimal prediction $\x_{k+1|k}$ can be bounded by using standard results for the projected gradient method, for the proof see for instance~\cite{Ryu2015,Taylor2017}.

\vskip2mm

\begin{proposition}\label{prop.intheo.2}
{\bf
Under Assumption~\ref{as.str}, we have that
\begin{equation}\label{init1}
\|\tilde{\x}_{k+1|k} - \x_{k+1|k}\| \leq \varrho_{\mathrm{P}}^{P}\,\|{\x}_{k} - \x_{k+1|k}\|, 
\end{equation}
with $\varrho_{\mathrm{P}} = \max\{|1-\alpha m|,|1 - \alpha L|\}$.
}
\end{proposition}

\vskip2mm
By putting together Propositions~\ref{prop.intheo.1} and~\ref{prop.intheo.2} and~\eqref{eq.lip} , we obtain for the total error after prediction as
\begin{subequations}
\begin{align}
\|\tilde{\x}_{k+1|k} \!- \!\x^*(t_{k+1})\| & \!\leq \!\|\tilde{\x}_{k+1|k} \!-\! \x_{k+1|k}\| \!+\! \|\x_{k+1|k}\!- \!\x^*(t_{k+1})\|\\ 
&\hspace*{-3.25cm}\leq \varrho_{\textrm{P}}^{P}\,\|{\x}_{k} \!-\! \x_{k+1|k}\| \!+\! \|\x_{k+1|k}\!-\! \x^*(t_{k+1})\|\\ 
&\hspace*{-3.25cm}\leq \varrho_{\textrm{P}}^{P}\,(\|{\x}_{k}\! - \!\x^*(t_k)\| \! +\!\|\x^*(t_k) \!-\!\x^*(t_{k+1})\| +\nonumber \\ 
& \hspace*{-2.25cm}\!\|\x^*(t_{k+1})\! - \!\x_{k+1|k}\|)\! +\! \|\x_{k+1|k}\!-\! \x^*(t_{k+1})\| \\
&\hspace*{-3.25cm}\leq \varrho_{\textrm{P}}^{P}\,\|{\x}_{k} - \x^*(t_k)\| + \Big(\varrho_{\textrm{P}}^{P}+1\Big) \|\x_{k+1|k}- \x^*(t_{k+1})\|  +\nonumber \\ 
& \hspace*{-2.25cm}\varrho_{\textrm{P}}^{P}\| \x^*(t_{k+1}) -  \x^*(t_{k})\|\\ \label{final.result}
&\hspace*{-3.25cm}\leq \eta_0 \|{\x}_{k} - \x^*(t_k)\|^2 +  \eta_1 \|{\x}_{k} - \x^*(t_k)\| + \eta_2,
\end{align}
\end{subequations}
where we have defined
\begin{align} \textrm{Case a)} &\quad \left\{
\begin{array}{rl}
\eta_0 &= 0, \\
\eta_1 &= \varrho_{\textrm{P}}^{P} + (\varrho_{\textrm{P}}^{P}+1) \frac{2L}{m},\\
\eta_2 &= (2 \varrho_{\textrm{P}}^{P} + 1) \Delta_1, 
\end{array}\right.
\\
\textrm{Case b)}& \quad \left\{
\begin{array}{rl}
\eta_0 &= (\varrho_{\textrm{P}}^{P}+1) \frac{C_1}{2\,m}, \\
\eta_1 &= \varrho_{\textrm{P}}^{P} + h\,(\varrho_{\textrm{P}}^{P}+1) \Big(\frac{C_1 C_0}{m^2} + \frac{C_2}{m}\Big),\\
\eta_2 &= \varrho_{\textrm{P}}^{P} \Big(h\,\frac{C_0}{m} + \Delta_2\Big) + \Delta_2, 
\end{array}\right.
\end{align}
where $\Delta_1$ and $\Delta_2$ are defined as in Proposition~\ref{prop.err}.
\vskip2mm

\noindent \textbf{Correction error.} We look now at the correction step, which by using standard results for the projected gradient method, we have
\begin{equation}\label{corr_eq}
\|{\x}_{k+1} - \x^*(t_{k+1})\| \leq \varrho_{\textrm{C}}^{C}\,\|\tilde{\x}_{k+1|k} - \x^*(t_{k+1})\|, 
\end{equation}
with $\varrho_{\textrm{C}} = \max\{|1-\alpha m|,|1 - \alpha L|\}$. And by putting together the result~\eqref{final.result} with~\eqref{corr_eq}, we obtain the recursive error bound,
\begin{multline}\label{rec}
\|{\x}_{k+1} - \x^*(t_{k+1})\| \leq \varrho_{\textrm{C}}^{C}\Big(\eta_0 \|{\x}_{k} - \x^*(t_k)\|^2 + \\ \eta_1 \|{\x}_{k} - \x^*(t_k)\| + \eta_2\Big).
\end{multline}

\vskip2mm

\noindent \textbf{Global error and convergence.} 

We start with Case a).
Call for simplicity $\bar{\eta}_i = \varrho_{\textrm{C}}^{C}\eta_i$, for $i = 0,1,2$. Since $\eta_0 = 0$, then convergence is achieved if 
\begin{equation}
\bar{\eta}_1 < 1 \iff \tau_0 = \varrho_{\textrm{C}}^{C}\left[\varrho_{\textrm{P}}^{P} + (\varrho_{\textrm{P}}^{P}+1) \frac{2L}{m}\right] < 1,
\end{equation}
which guarantees a monotonical decrease of the optimality gap with linear rate and asymptotic bound of
\begin{equation}
\limsup_{k \to \infty}\|{\x}_{k+1} - \x^*(t_{k+1})\|  = \frac{\varrho_{\textrm{C}}^{C}\Delta_1}{1-\tau_0} (2 \varrho_{\textrm{P}}^{P} + 1) = O(\varrho_{\textrm{C}}^{C} h),
\end{equation}
which are the claims in Theorem~\ref{theo.constrained_0}.

\vskip2mm
We finish with Case b).
Call again for simplicity $\bar{\eta}_i = \varrho_{\textrm{C}}^{C}\eta_i$, for $i = 0,1,2$. Then convergence is achieved if
\begin{enumerate}
\item[1)] Each iteration does not increase the error, so that
\begin{multline}\label{decrease}
\bar{\eta}_0 \|{\x}_{k} - \x^*(t_k)\|^2 +  \bar{\eta}_1 \|{\x}_{k} - \x^*(t_k)\| + \bar{\eta}_2 \leq \\ \tau \|{\x}_{k} - \x^*(t_k)\| + \bar{\eta}_2 
\end{multline}
for a $\tau < 1$; 
\item[2)] One can find a $\tau <1$ such that~\eqref{decrease} holds. 
\end{enumerate}
By simple algebra, convergence is achieved if $\bar{\eta}_1 < \tau < 1$, that is if 
\begin{equation}
\varrho_{\textrm{P}}^{P}\varrho_{\textrm{C}}^{C} < \tau, \quad h < \frac{\tau - \varrho_{\textrm{C}}^{C}\varrho_{\textrm{P}}^{P}}{\varrho_{\textrm{C}}^{C}(\varrho_{\textrm{P}}^{P} + 1)}\Big(\frac{C_1 C_0}{m^2} + \frac{C_2}{m} \Big)^{-1} = \bar{h},
\end{equation}
which sets the bounds on the number of prediction and correction steps as well as the sampling period, and if the initial optimality gap is chosen as
\begin{equation}
\|\x_0 - \x^*(t_0) \| \leq \frac{\tau - \bar{\eta}_1}{\bar{\eta}_0} = \bar{R}.
\end{equation}
The convergence region depends on the sampling period and on the prediction and correction steps. When $h \to 0$, then
\begin{equation}
\lim_{h\to 0} \bar{R} = \frac{2\, m}{C_1}\,\frac{\tau - \varrho_{\textrm{P}}^{P}\varrho_{\textrm{C}}^{C}}{\varrho_{\textrm{C}}^{C}(\varrho_{\textrm{P}}^{P} + 1)}. 
\end{equation}

As for the convergence asymptotical error, by using~\eqref{rec} in combination with~\eqref{decrease}, we can show that
\begin{equation}\label{final}
\|\x_k - \x^*(t_k)\| \leq \tau^{k} \|\x_0 - \x^*(t_0)\| + \bar{\eta}_2 \frac{1-\tau^{k}}{1 - \tau},
\end{equation}
from which, by letting $k \to \infty$, the result~\eqref{main.result1} follows. \qed


\section{Proof of Theorem~\ref{theo.unconstrained} }\label{ap.newton}

We proceed similarly to Appendix~\ref{ap.grad}. One of the main difference is a new Proposition~\ref{prop.intheo.1} suited for the situation at hand. 

\vskip2mm 

\begin{proposition}\label{prop.intheo.1bis}
{\bf
The following facts hold true. 
\begin{enumerate}
\item[{Case a)}] Under the same assumptions and notation of Theorem~\ref{theo.unconstrained_0}, we have that
\begin{equation}\label{pp_0}
\hspace*{-.5cm}\|\x_{k+1|k}\! - \!\x^*(t_{k+1})\| \!\leq\! \Big(1-\gamma + \gamma \frac{2 L}{m}\Big) \|\x_k\!-\! \x^*(t_k)\|\!+\! 2\Delta_1.
\end{equation}
\item[{Case b)}] Under the same assumptions and notation of Theorem~\ref{theo.unconstrained}, we have that
\begin{multline}\label{pp}
\|\x_{k+1|k} - \x^*(t_{k+1})\| \leq \gamma\,\frac{C_1}{2\,m}\|\x_k- \x^*(t_k)\|^2 + \\\Big[1-\gamma + h\,\Big(\frac{C_1 C_0}{m^2} + \frac{C_2}{m}  \Big)\Big]\|\x_k- \x^*(t_k)\|+ \Delta_2.
\end{multline}
\end{enumerate}
}
\end{proposition}
\vskip2mm

\begin{IEEEproof}
Start by noticing that Proposition~\ref{prop.err} holds true even when $\gamma < 1$ (in the unconstrained case), since $\nabla_{\x}f(\x^*(t_k); t_k) = \nabla_{\x}f_0 = \textbf{0}$. In fact, Eq.~\eqref{proofA.dummy1} should read
\begin{equation}
\nabla_{\x} f_0 + \Q_0 (\x - \x_0) + h\, \c_0 = (1-\gamma)\nabla_{\x} f_0,
\end{equation} 
but this is in fact equivalent to the original \eqref{proofA.dummy1}, since $\nabla_{\x}f_0 = \textbf{0}$, and therefore the whole proposition is still valid, and in particular $\|\x^*_{k+1|k} - \x^*(t_{k+1})\| \leq \Delta_i$, where $i=1$ for Case a) and $i=2$ for Case b).

Call now $\delta\x_k = \x_{k+1|k} - \x_k$, and $\delta\x^*_k = \x^*_{k+1|k} - \x^*(t_k)$. Then, 
\begin{multline}
\|\x_{k+1|k} - \x^*(t_{k+1})\| = \\ \|\x_{k} + \delta\x_k - (\x^*(t_k) + \delta\x^*_k) + ( \x^*_{k+1|k} - \x^*(t_{k+1})) \|,
\end{multline}
which given Proposition~\ref{prop.err} and by using the Triangle inequality can be upper bounded as
\begin{equation}\label{b_1}
\|\x_{k+1|k} \!-\! \x^*(t_{k+1})\|\! \leq \!(1-\gamma)\|\x_{k}\! -\! \x^*(t_k)\| + \|\delta\check{\x}_k \!-\!  \delta\x^*_k\| \!+\! \Delta_i,
\end{equation}
where we have set $\delta\check{\x}_k = \delta\x_k + \gamma(\x_k - \x^*(t_k))$. Note that this decomposition may seem cumbersome, yet it is the cornerstone of the proof of this proposition.  

For $\delta\check{\x}_k$ and $\delta \x^*_k$, it holds that
\begin{align}
& \gamma\, \nabla_{\x} f_k + \gamma\, \Q_k (\x_0 - \x_k) + \Q_k \delta\check{\x}_k + h\, \c_k = \mathbf{0},\\
& \gamma \nabla_{\x} f_0 + \Q_0\delta{\x}^*_k + h\, \c_0 = \mathbf{0},
\end{align}
where as in~\eqref{simplif2} we have used the simplifications
\begin{subequations}
\begin{align}
\nabla_{\x} f_k = \nabla_{\x} f (\x_k; t_{k}), \quad & \Q_k = \nabla_{\x\x} f (\x_k; t_{k}) \\
\c_k = \nabla_{t\x} f (\x_k; t_{k}),\quad &  \x = \x_{k+1|k}.
\end{align}
\end{subequations}
while $\nabla_{\x} f_0$, $\Q_0$, $\c_0$, and $\x_0$ are defined just as in~\eqref{simplif}-\eqref{simplif22}, for $i = 0$ (i.e., $\x_0 = \x^*(t_k)$, and so on). Define
\begin{align}
g(\delta\x) &= \gamma\, (\nabla_{\x} f_k - \nabla_{\x} f_0 + \Q_k (\x_0 - \x_k)) + \nonumber \\& \hspace*{1cm} (\Q_k - \Q_0)\delta\x + h\, (\c_k-\c_0), \\
G(\delta\x) &= \gamma \nabla_{\x} f_0 + \Q_0\delta{\x} + h\, \c_0,
\end{align}
and notice that $(g + G)(\delta\check{\x}_k) = \mathbf{0}$, while $(g + G)(\delta{\x}^*_k) = g(\delta{\x}^*_k)$. With a similar argument as the one of the proof of Proposition~\ref{prop.intheo.1}, then
\begin{equation}\label{b_0}
\|\delta\check{\x}_k - \delta{\x}^*_k\| \leq \frac{1}{m} \|g(\delta{\x}^*_k)\|.  
\end{equation}
Let us now bound $\|g(\delta{\x}^*_k)\|$. 

Case a) Due to Assumptions~\ref{as.str}:
\begin{equation}\label{b_2_0}
\|g(\delta{\x}^*_k)\| \!\leq \! 2 \gamma L \|\x_k \!-\! \x_0\| + 2L\|\delta{\x}^*_k\| + 2 h\, C_0.
\end{equation}  

Case b) Due to Assumptions~\ref{as.str}-\ref{as.smooth}, 
\begin{equation}\label{b_2}
\|g(\delta{\x}^*_k)\| \!\leq \!\gamma \frac{C_1}{2}\|\x_k \!-\! \x_0\|^2 + C_1  \|\x_k - \x_0\|\|\delta{\x}^*_k\| + h\, C_2 \|\x_k \!-\! \x_0\|.  
\end{equation}

The next step of the proof is to upper bound the term $\|\delta \x^*_k\| = \|\x^*_{k+1|k} - \x^*(t_k)\|$. We know that
\begin{equation}
\gamma \nabla_{\x}f_0 + \Q_0\delta \x^*_k + h \c_0 = \mathbf{0},
\end{equation}
and since $\nabla_{\x}f_0  = \mathbf{0}$, then $\delta \x^*_k = - h\,\Q_0^{-1} \c_0$. Which yields 
\begin{equation}\label{b_3}
\|\delta \x^*_k\| \leq h\,\frac{C_0}{m}.
\end{equation}
By putting together the bounds~\eqref{b_1}, \eqref{b_0}, \eqref{b_2_0}, \eqref{b_2}, and \eqref{b_3}, the claims~\eqref{pp_0} and \eqref{pp} follow. 
\end{IEEEproof}

\vskip2mm
By using Proposition~\ref{prop.intheo.1bis} along with the same arguments as the one in Eq.s~\eqref{init1} till \eqref{rec}, we arrive at the error recursion
\begin{equation}
\|{\x}_{k+1} - \x^*(t_{k+1})\| \leq \bar{\eta}_0' \|{\x}_{k} - \x^*(t_k)\|^2 +  \bar{\eta}_1' \|{\x}_{k} - \x^*(t_k)\| +  \bar{\eta}_2',
\end{equation}
where,
\begin{equation}
\textrm{Case a)}\, \left\{
\begin{array}{rl}
\bar{\eta}_0' &= 0, \\
\bar{\eta}_1' &= \varrho_{\textrm{C}}^{C} \left[\varrho_{\textrm{P}}^{P} +  (\varrho_{\textrm{P}}^{P}+1)\Big(1-\gamma +\gamma \frac{2L}{m}\Big)\right],\\
\bar{\eta}_2' &=  \varrho_{\textrm{C}}^{C} \left[2 (2 \varrho_{\textrm{P}}^{P} + 1) \Delta_1\right],
\end{array}\right.
\end{equation}
\begin{multline}
\hspace*{-1mm}\textrm{Case b)}\\ \left\{
\begin{array}{rl}
\bar{\eta}_0' &= \gamma \varrho_{\textrm{C}}^{C}(\varrho_{\textrm{P}}^{P}+1) \frac{C_1}{2\,m}, \\
\bar{\eta}_1' &= \varrho_{\textrm{C}}^{C}\Big[\varrho_{\textrm{P}}^{P} + (\varrho_{\textrm{P}}^{P}+1)\Big( 1-\gamma + h\, \Big(\frac{C_1 C_0}{m^2} + \frac{C_2}{m}\Big)\Big)\Big],\\
\bar{\eta}_2' &= \varrho_{\textrm{C}}^{C}\Big[\varrho_{\textrm{P}}^{P} \Big(h\,\frac{C_0}{m} + \Delta_2\Big) + \Delta_2\Big]. 
\end{array}\right.
\end{multline}
\vskip2mm
By using now the same reasoning as in Eq.s~\eqref{decrease} till Eq.s~\eqref{final}, the claims~\eqref{main.result2_0}-\eqref{main.result2} can be proven.  
\qed


\bibliographystyle{ieeetr}

\bibliography{PaperCollection_0}

\end{document}